\documentclass[12pt]{article}
\usepackage{latexsym,amsfonts,amssymb}
\usepackage{leftindex}
\usepackage{bm}
\usepackage[dvips]{color}
\usepackage{colordvi,multicol}
\usepackage{amsmath}
 \usepackage{graphicx}
 \usepackage[english]{babel}
\usepackage{caption}
\usepackage{tabu}
\usepackage{float}
\usepackage[pagewise, left, displaymath, mathlines]{lineno}
\usepackage{dsfont}
\usepackage{xcolor}
\usepackage{tikz-cd}
\usepackage[shortlabels]{enumitem}

\newcommand{\xms}{X_{t-}^\star}
\newcommand{\yms}{Y_{t-}^\star}
\newcommand{\zms}{Z_{t-}^\star}
\newcommand{\xd}{X_{t}^\dagger}
\newcommand{\yd}{Y_{t}^\dagger}
\newcommand{\zd}{Z_{t}^\dagger}

\def \cal{\mathcal}

\newtheorem{thm}{Theorem}[section]

\newtheorem{rem}[thm]{Remark}
\newtheorem{exa}[thm]{Example}
\newtheorem{simu}[thm]{Simulation}

\date{}

\begin{document}

\title{\bf A Unified Stochastic SIR Model Driven By L\'{e}vy Noise With Time-Dependency}
\author{Terry Easlick  and Wei Sun\thanks{Corresponding author.}\\ \\
  {\small Centre de recherche du CHU Sainte-Justine}\\
{\small D\'epartement de Math\'ematiques et de Statistique}\\
{\small Universit\'e de Montr\'eal, Canada}\\ \\
 {\small Department of Mathematics and Statistics}\\
{\small Concordia University, Canada}\\ \\
{\small  terry.easlick@umontreal.ca,\ \ \ \ wei.sun@concordia.ca}}

\maketitle

\begin{abstract}
\noindent We propose a unified stochastic SIR model driven by L\'{e}vy noise. The model is structural enough to allow
for time-dependency, nonlinearity, discontinuity, demography and environmental disturbances. We present concise results on the existence and uniqueness of positive global
solutions and investigate the extinction and persistence of the novel model. Examples and simulations are provided to illustrate the main results.
\end{abstract}



\section{Introduction}

The investigation of infectious diseases has a rich history that gave rise to the field epidemiology (cf. \cite{ephist}). The first known study of infectious disease data relating to causes of death occurred in the 17th century by John Graunt. The work completed by Graunt was published in his 1662 book {\it Natural and Political Observations made upon the Bills of Mortality}.  Mathematical epidemiology is considered to have begun with the work of Daniel Bernoulli with his model for smallpox inoculation in the 18th century. Late in the early 20th century, W.H. Hamer surmised the rate of disease spread was dependent upon the numbers of those at risk (susceptible) and those currently ill/contagious (infected). He also recommended the use of a mass action law for the rate of infections. The two previous ideas cemented the basic building blocks of epidemiological compartmental models. Compartmental models are the idea that the population can be segmented based on the status of health as related to a disease. For instance, the subpopulation at risk of becoming infected by a particular disease is represented by the susceptible compartment where as the subpopulation which is currently infected is represented by the infected compartment.

Compartmental models have garnered much attention in the past century by researchers in pursuit of understanding and control of infectious diseases. Mathematical analysis of such models aids decision-making regarding public health policy changes -- especially in the event of a pandemic (e.g., COVID-19).  One such compartmental model introduced by Kermack and McKendrick \cite{KM} in 1927 divides a population into three compartments -- susceptible, infected, and recovered (SIR).   The classical SIR model is as follows:
\begin{eqnarray}\label{SIRA}
\begin{cases}
	\frac{dX_t}{dt} =  - \beta X_tY_t, \\
	\frac{dY_t}{dt} = \left( \beta X_t - \gamma\right)Y_t,\\
	\frac{dZ_t}{dt} = \gamma Y_t,
\end{cases}
\end{eqnarray}
where $\beta$ is the transmission coefficient and $\gamma$ the recovery coefficient. Additionally, demography may be introduced to include birth recruitment coefficient $\Lambda$ and mortality coefficient $\mu$ as:
\begin{eqnarray}\label{SIRB}
\begin{cases}
	\frac{dX_t}{dt} = \Lambda - \mu X_t - \beta X_tY_t, \\
	\frac{dY_t}{dt} = \left[ \beta X_t - (\mu+\gamma )\right]Y_t,\\
	\frac{dZ_t}{dt} = \gamma Y_t - \mu Z_t.
\end{cases}
\end{eqnarray}

The basic SIR models (\ref{SIRA}) and (\ref{SIRB}) have many variations including the SIRD, SIRS, SIRV, SEIR, MSIR, etc. (cf. e.g., \cite{Ba}, \cite{BR} and \cite{SK}). Further, these deterministic models have been put into different
stochastic frameworks, which make the situation more realistic (cf. e.g., \cite{Bao},
\cite{Cara}--\cite{G}, \cite{JJS}, \cite{Liu}, \cite{Pri} and \cite{T}--\cite{Zhou2}). The existing models are often analyzed with a focus on specific diseases or parameters. Such studies have been very successful in achieving new results; however, often it is the case that structural variability is lacking in these models. To overcome the drawbacks inherent in traditional approaches, we propose and investigate in this paper the unified stochastic SIR (USSIR) model:
\begin{eqnarray}\label{p1_eq2}
\begin{cases}
\begin{split}
&dX_t =\ b_1(t,X_t,Y_t,Z_t)dt + \sum\limits_{j=1}^n \sigma_{1j}(t,X_t,Y_t,Z_t)dB_t^{(j)}\\
	&\indent\indent + \int_{\left\{\lvert u \rvert < 1\right\}}\hspace{-.3cm}h_1(t,X_{t-},Y_{t-},Z_{t-},u)\tilde{N}(dt,du)+ \int_{\left\{\lvert u \rvert \ge 1\right\}}\hspace{-.3cm}g_1(t,X_{t-},Y_{t-},Z_{t-},u)N(dt,du),
\end{split}\\
\begin{split}
&dY_t =\ b_2(t,X_t,Y_t,Z_t)dt + \sum\limits_{j=1}^n\sigma_{2j}(t,X_t,Y_t,Z_t)dB_t^{(j)}\\
	&\indent\indent+ \int_{\left\{\lvert u \rvert < 1\right\}}\hspace{-.3cm}h_2(t,X_{t-},Y_{t-},Z_{t-},u)\tilde{N}(dt,du)+  \int_{\left\{\lvert u \rvert \ge 1\right\}}\hspace{-.3cm}g_2(t,X_{t-},Y_{t-},Z_{t-},u)N(dt,du),
\end{split}\\
\begin{split}
&dZ_t =\ b_3(t,X_t,Y_t,Z_t)dt + \sum\limits_{j=1}^n\sigma_{3j}(t,X_t,Y_t,Z_t)dB_t^{(j)} \\
	&\indent\indent+\int_{\left\{\lvert u \rvert < 1\right\}}\hspace{-.3cm}h_3(t,X_{t-},Y_{t-},Z_{t-},u)\tilde{N}(dt,du) + \int_{\left\{\lvert u \rvert \ge 1\right\}}\hspace{-.3cm}g_3(t,X_{t-},Y_{t-},Z_{t-},u)N(dt,du).
\end{split}\\
\end{cases}
\end{eqnarray}
Hereafter, $\mathbb{R}_+$ denotes the set of all positive real numbers, $\left(B^{(1)}_t,\dots,B^{(n)}_t\right)_{t \ge 0}$ is a standard $n$-dimensional Brownian motion, $N$ is a Poisson random measure on
$\mathbb{R}_+\times (\mathbb{R}^d-\{0\})$ with intensity measure $\nu$ satisfying $\int_{\mathbb{R}^d-\{0\}}(1\wedge|u|^2)\nu(du)<\infty$ and $\tilde{N}(dt,du) = N(dt,du) - \nu(du)dt$, $\left(B^{(1)}_t,\dots,B^{(n)}_t\right)_{t \ge 0}$ and $N$ are independent, $b_i, \sigma_{ij}:[0,\infty) \times \mathbb{R}^3_+\mapsto\mathbb{R}$, $h_i,g_i:  [0,\infty) \times \mathbb{R}^3_+\times (\mathbb{R}^d-\{0\}) \mapsto\mathbb{R}$, $i=1,2,3,\ j=1,2,\ldots,n$, are measurable functions.

 We will show that the USSIR  model (\ref{p1_eq2}) is structural in design that allows variability without sacrificing key results on the extinction and persistence of diseases. Namely, the model allows for time-dependency, nonlinearity (of drift, diffusion and jump) and demography. Environmental disturbances can have profound effects on transmission, recovery, mortality and population growth. Different diseases may have drastically different dynamics thus making perturbations of transmission or recovery natural to consider. Moreover, exact mixing of susceptible and infected compartments may be difficult to measure thus the transmission will be prone to disturbances -- this can be accounted for by the inclusion of stochasticity. The above model encapsulates the stochastic perturbations driven by white noises $(B^{(j)}_t)$ with intensities $\sigma_{ij}(t,X_t,Y_t,Z_t)$ and Poisson random measure $N(dt,du)$ with small jumps $h_{i}(t,X_{t-},Y_{t-},Z_{t-},u)$ and large jumps $g_{i}(t,X_{t-},Y_{t-},Z_{t-},u)$. An important structural feature we emphasize is time-dependency. Time-dependency can capture the progression of a disease  insofar as mutations/transmissibility (e.g., Delta and Omicron variants of COVID-19, vaccination programs).

The remainder of this paper is organized as follows. In Sections 2 and 3, we establish results on the existence and uniqueness of global solutions, extinction and persistence of diseases, and provide illustrative examples. Section 2  is concerned with the USSIR model (\ref{p1_eq2}) where compartments take values in $\mathbb{R}^3_+$ and  Section 3  is concerned with a special case  focusing on proportionality. The examples in Sections 2 and 3 are intended to show the flexibility of the model while maintaining biological relevance. In section 4, we present simulations which correspond to examples given in the previous sections. In section 5, we make concluding remarks. At the time of writing this paper we are unaware of existing work on the USSIR model and aim to add such a novel model to the existing literature.

\section{Model for population numbers}\setcounter{equation}{0}
In this section, we let $X_t$, $Y_t$ and $Z_t$ denote respectively the numbers of susceptible, infected and recovered individuals at time $t$. For $t\in[0, \infty)$, $(x,y,z)\in\mathbb{R}^3_+$ and $u\in \mathbb{R}^d-\{0\}$, define
\begin{eqnarray}\label{p1_eq3}
&&k(t,x,y,z,u):= \frac{h_1(t,x,y,z,u)}{x}+\frac{h_2(t,x,y,z,u)}{y}+\frac{h_3(t,x,y,z,u)}{z}\nonumber\\
&&\ \ \ \ \ \ - \ln\left\{\left( 1+ \frac{h_1(t,x,y,z,u)}{x}\right)\left( 1+ \frac{h_2(t,x,y,z,u)}{y}\right)\left( 1+ \frac{h_3(t,x,y,z,u)}{z}\right)\right\}.\ \ \ \ \ \ \ \
\end{eqnarray}
We make the following assumptions.

\begin{enumerate}[\bf({A}1)]

\item\label{assm21} There exists $(x_0,y_0,z_0)\in\mathbb{R}^3_+$ such that for any $T \in (0,\infty)$ and $i=1,2,3$,
\begin{flalign*}
	&b_i(\cdot,x_0,y_0,z_0),\ \ \sum\limits_{j=1}^n|\sigma_{ij}(\cdot,x_0,y_0,z_0)| \in L^2[0,T],\ \ \int_{\{|u|<1\}} |h_i(\cdot,x_0,y_0,z_0,u)|^2\nu(du) \in L^1[0,T].
\end{flalign*}

\item\label{assm22} For any $T \in (0,\infty)$ and $N \in \mathbb{N}$, there exists $K_{N,T} \in L_+^1[0,T]$ such that
\begin{eqnarray*}
&&\sum\limits_{i=1}^3  \lvert b_i(t,x_1,y_1,z_1)-b_i(t,x_2,y_2,z_2)\rvert^2 +\sum\limits_{i=1}^3\sum\limits_{j=1}^n\lvert \sigma_{ij}(t,x_1,y_1,z_1)-\sigma_{ij}(t,x_2,y_2,z_2)\rvert^2\\
	&&+\sum\limits_{i=1}^3\int_{\{|u|<1\}}| h_i(t,x_1, y_1,z_1,u) - h_i(t,x_2,y_2,z_2,u)|^2\nu(du)\\
& \le& K_{N,T}(t)[(x_1-x_2)^2+(y_1-y_2)^2+(z_1-z_2)^2],\\
&&\forall t\in [0,T],\ (x_1,y_1,z_1),(x_2,y_2,z_2)\in \left[\frac{1}{N},N\right]^3.
\end{eqnarray*}

\item\label{assm23} For any $(x,y,z) \in \mathbb{R}^3_+$, $t\in(0,\infty)$ and $u\in \mathbb{R}^d-\{0\}$,
\begin{eqnarray*}
&&\left( 1+ \frac{h_1(t,x,y,z,u)}{x}\right),\ \left( 1+ \frac{h_2(t,x,y,z,u)}{y}\right),\ \left( 1+ \frac{h_3(t,x,y,z,u)}{z}\right) > 0,
\end{eqnarray*}
and
\begin{eqnarray*}
&&\left( 1+ \frac{g_1(t,x,y,z,u)}{x}\right),\ \left( 1+ \frac{g_2(t,x,y,z,u)}{y}\right),\ \left( 1+ \frac{g_3(t,x,y,z,u)}{z}\right) > 0.
\end{eqnarray*}

\item\label{assm24} For any $T\in (0,\infty)$,
\begin{eqnarray*}
&&\sup_{t\in[0, T],\, (x,y,z)\in\mathbb{R}^3_+}\left\{\frac{(x-1)b_1(t,x,y,z)}{x}+\frac{(y-1)b_2(t,x,y,z)}{y}+\frac{(z-1)b_3(t,x,y,z)}{z}\right\}<\infty,
\end{eqnarray*}
\begin{eqnarray*}
&&\sum_{j=1}^n\sup_{t\in[0, T],\, (x,y,z)\in\mathbb{R}^3_+}\left\{\frac{|\sigma_{1j}(t,x,y,z)|}{x}+\frac{|\sigma_{2j}(t,x,y,z)|}{y}+\frac{|\sigma_{3j}(t,x,y,z)|}{z}\right\}<\infty,
\end{eqnarray*}
and there exists $\eta_T\in {\cal B}(\mathbb{R}^d)$ such that
\begin{eqnarray}\label{p1_eq4}
\eta_T>0,\ \ \ \ \int_{\{|u|<1\}}\eta_T(u)\nu(du)<\infty,
\end{eqnarray}
and
$$
\sup_{t\in[0, T],\, (x,y,z)\in\mathbb{R}^3_+, u\in \mathbb{R}^d-\{0\}}\frac{k(t,x,y,z,u)}{\eta_T(u)}<\infty.
$$

\end{enumerate}

Let $C^{1,2}(\mathbb{R}_+\times\mathbb{R}_+^3;\mathbb{R})$ be the
set of all functions $V(t,w)$ on
$\mathbb{R}_+\times\mathbb{R}_+^3$ which are continuously
differentiable with respect to $t$ and twice continuously
differentiable with respect to $w=(w_1,w_2,w_2)$. For $V\in
C^{1,2}(\mathbb{R}_+\times\mathbb{R}_+^3;\mathbb{R})$, we define
\begin{eqnarray*}
{\cal L}V(t,w)&=&\langle V_w(w), b(t,w)\rangle+
\frac{1}{2}\mathrm{ trace}(\sigma^T(t,w)V_{ww}(w)\sigma(t,w))\nonumber\\
&&+\int_{\{|u|<1\}}[V(w+h(t,w,u))-V(w)-
\langle V_w(w), h(t,w,u)\rangle]\nu(\mathrm{d}u),
\end{eqnarray*}
with $V_w=(\frac{\partial V}{\partial w_1},\frac{\partial V}{\partial w_2},\frac{\partial V}{\partial w_3})$,
$V_{ww}=(\frac{\partial^2 V}{\partial w_i\partial w_j})_{1\le i,j\le 3}$, $b=(b_1,b_2,b_3)$, $\sigma=(\sigma_{ij})_{1\le i,j\le 3}$ and $h=(h_1,h_2,h_3)$.

Now we present the result on the existence and uniqueness of solutions to the system (\ref{p1_eq2}).
\begin{thm}\label{RE}
Suppose that Assumptions \ref{assm21}--\ref{assm24} hold. Then, for any given initial value $(X_0,Y_0, Z_0) \in \mathbb{R}^3_+$, the system (\ref{p1_eq2}) has a unique strong solution taking values in $\mathbb{R}^3_+$.
\end{thm}

\noindent {\bf Proof.}\ \ By \ref{assm23} and the interlacing technique (cf. \cite{A}), to complete the proof, we need only  consider the case that $g_i\equiv0$, $i=1,2,3$. Then, equation (\ref{p1_eq2}) becomes equation
\begin{eqnarray}\label{p1_eq5}
\begin{cases}

dX_t =\ b_1(t,X_t,Y_t,Z_t) + \sum\limits_{j=1}^n \sigma_{1j}(t,X_t,Y_t,Z_t)dB_t^{(j)} + \int_{\left\{\lvert u \rvert < 1\right\}}h_1(t,X_{t-},Y_{t-},Z_{t-},u)\tilde{N}(dt,du),\\
dY_t =\ b_2(t,X_t,Y_t,Z_t) + \sum\limits_{j=1}^n\sigma_{2j}(t,X_t,Y_t,Z_t)dB_t^{(j)} + \int_{\left\{\lvert u \rvert < 1\right\}}h_2(t,X_{t-},Y_{t-},Z_{t-},u)\tilde{N}(dt,du),\\
dZ_t =\ b_3(t,X_t,Y_t,Z_t)+ \sum\limits_{j=1}^n\sigma_{3j}(t,X_t,Y_t,Z_t)dB_t^{(j)} +\int_{\left\{\lvert u \rvert < 1\right\}}h_3(t,X_{t-},Y_{t-},Z_{t-},u)\tilde{N}(dt,du).
\end{cases}
\end{eqnarray}

 By \ref{assm21} and \ref{assm22}, similar to \cite[Lemma 2.1]{Sun}, we can show that there exists a unique local strong solution to equation (\ref{p1_eq5}) on $[0,\tau)$, where $\tau$ is the explosion time. We will show below that $\tau = \infty$ almost surely (a.s.). Define
\begin{eqnarray*}
\tau_N =\inf\left\{t\in[0,\tau): (X_t,Y_t,Z_t)\not\in\left[\frac{1}{N},N\right]^3\right\},\ \ \ \ N\in\mathbb{N},
\end{eqnarray*}
and
$$
\tau_{\infty}=\lim_{N\rightarrow\infty}\tau_N.
$$
We have that $\tau_{\infty} \le \tau$ so it suffices to show $\tau_{\infty} = \infty$ a.s.. Hence assume the contrary that there exist $\varepsilon>0$ and $T>0$ such that
$$
\mathbb{P}(\tau_{\infty} < T) >\varepsilon,
$$
which implies that
\begin{eqnarray}\label{C11}
\mathbb{P}(\tau_N < T) >\varepsilon,\ \ \ \ \forall N\in\mathbb{N}.
\end{eqnarray}

Define
$$
V(x,y,z) =(x-1-\ln x)+ (y-1-\ln y)+(z-1-\ln z),\ \ \ \ (x,y,z)\in (0,\infty)^3,
$$
and
$$
W_t=(X_t,Y_t,Z_t).
$$
By It\^{o}'s formula, we obtain that for $t\le \tau_N$,
\begin{eqnarray*}
V(W_t)&=& V(W_0)+\int_0^t\mathcal{L}V(s,W_s)ds +\int_0^t \langle V_x(W_s), \sigma(s,X_s) \rangle dB_s\\
&&+ \int_0^t\int_{\{|u|<1\}}\left[V(W_{s-}+h(s,W_{s-},u)) - V(W_{s-}) \right]\tilde{N}(ds,du).
\end{eqnarray*}
 Then, by \ref{assm24}, there exist $C_T>0$ and  $\eta_T\in {\cal B}(\mathbb{R}^d)$ such that (\ref{p1_eq4}) holds and
\begin{eqnarray*}
&&\mathbb{E}\left[V(X_{T\wedge \tau_N},Y_{T\wedge \tau_N},Z_{T\wedge \tau_N})\right] -V(X_0,Y_0,Z_0)\\
&=& \mathbb{E}\left[\int_0^{T\wedge \tau_N} \left\{\frac{(X_s-1)b_1(s,X_s,Y_s,Z_s)}{X_s}+\frac{(Y_s-1)b_2(s,X_s,Y_s,Z_s)}{Y_s}+\frac{(Z_s-1)b_3(s,X_s,Y_s,Z_s)}{Z_s} \right\}ds\right]\\
	&&+\ \frac{1}{2}\sum\limits_{j=1}^n\mathbb{E}\left[\int_0^{T\wedge \tau_N} \left\{\frac{\sigma_{1j}^2(s,X_s,Y_s,Z_s)}{X_s^2}+\frac{\sigma_{2j}^2(s,X_s,Y_s,Z_s)}{Y_s^2}+\frac{\sigma_{3j}^2(s,X_s,Y_s,Z_s)}{Z_s^2}\right\}ds\right]\\
	&&+\mathbb{E}\left[\int_0^{T\wedge \tau_N}\ \int_{\{|u|<1\}}k(s,X_s,Y_s,Z_s,u)\nu(du)ds\right]\\
&\le& C_TT+\frac{C_T^2T}{2}+C_TT \int_{\{|u|<1\}}\eta_T(u)\nu(du).
\end{eqnarray*}
However, by (\ref{C11}), we get
$$
\mathbb{E}[V(X_{T\wedge\tau_N},Y_{T\wedge\tau_N},Z_{T\wedge\tau_N})]> \varepsilon\left[\left(\frac{1}{N}-1+\ln N\right)\wedge(N-1-\ln N)\right]\rightarrow\infty\ \ {\rm as}\ N\rightarrow\infty.
$$
We have arrived at a contradiction. Therefore, $\tau = \infty$ a.s. and the proof is complete.\hfill\fbox

Next, we consider the extinction and persistence of diseases. Namely, we investigate whether a disease will extinct with an exponential rate or will be persistent in mean. The system (\ref{p1_eq2}) is called persistent in mean if
\[
	\liminf_{t \to \infty} \frac{1}{t}\int_0^tY_sds > 0\ \ a.s..
\]

\begin{thm}\label{p1_thm2}

Suppose that Assumptions \ref{assm21}--\ref{assm24} hold. Let $(X_t,Y_t,Z_t)$ be a solution to equation (\ref{p1_eq2}) with $(X_0,Y_0,Z_0) \in\mathbb{R}^3_+$. We assume that
\begin{equation}\label{p1_eqa}
\int_0^{\infty}\frac{\varphi(t)}{(1+t)^2}dt<\infty,
\end{equation}
where
\begin{eqnarray*}
\varphi(t)&:=&\sup_{(x,y,z) \in \mathbb{R}^3_+} \left\{\frac{ \sum_{j=1}^n \sigma^2_{2j}(t,x,y,z)}{y^2}+\int_{\{|u|<1\}}\left[\ln\left(1+\frac{h_2(t,x,y,z,u)}{y}\right)\right]^2\nu(du)\right.\\
&&\ \ \ \ \ \ \ \ \ \ \ \left.+\int_{\{|u|\ge1\}}\left[\ln\left(1+\frac{g_2(t,x,y,z,u)}{y}\right)\right]^2\nu(du)\right\}.
\end{eqnarray*}

\noindent (i) If
\begin{eqnarray}\label{extqg}
\alpha&:=&\limsup\limits_{t\to\infty}\left\{\sup\limits_{(x,y,z)\in\mathbb{R}^3_+}\left[\frac{b_2(t,x,y,z)}{y} - \frac{\sum_{j=1}^n\sigma^2_{2j}(t,x,y,z)}{2y^2}\right]\right.\nonumber\\
&&\ \ \ \ \ \ \ \ \ \ \ \
+\int_{\{|u|<1\}}\sup\limits_{(x,y,z)\in\mathbb{R}^3_+}\left[\ln\left(1+\frac{h_2(t,x,y,z,u)}{y}\right)-\frac{h_2(t,x,y,z,u)}{y}\right]\nu(du)\nonumber\\
&&\ \ \ \ \ \ \ \ \ \ \ \
+\int_{\{|u|\ge1\}}\sup\limits_{(x,y,z)\in\mathbb{R}^3_+}\left[\ln\left(1+\frac{g_2(t,x,y,z,u)}{y}\right)\right]\nu(du)\nonumber\\
&<& 0,
\end{eqnarray}
then
\begin{eqnarray}\label{p1_eq9a}
\limsup_{t\to\infty} \frac{\ln Y_t}{t} \le \alpha\ a.s..
\end{eqnarray}

\noindent (ii) If there exist positive constants $\lambda_0$ and $\lambda$ such that
\begin{eqnarray}\label{per1a}
&&\liminf_{t\to\infty}\frac{1}{t}\int_0^t\bigg\{\lambda_0Y_s+\frac{b_2(s,X_s,Y_s,Z_s)}{Y_s} - \frac{\sum_{j=1}^n\sigma^2_{2j}(s,X_s,Y_s,Z_s)}{2Y_s^2}\nonumber\\
	&&\ \ \ \ \ \ \ \ \ \ \ \ \ \ +\int_{\{|u|<1\}}\left[\ln\left(1+\frac{h_2(s,X_{s-},Y_{s-},Z_{s-},u)}{Y_{s-}}\right)-\frac{h_2(s,X_{s-},Y_{s-},Z_{s-},u)}{Y_{s-}}\right]\nu(du)\nonumber\\
&&\ \ \ \ \ \ \ \ \ \ \ \ \ \ +\int_{\{|u|\ge1\}}\ln\left(1+\frac{g_2(s,X_{s-},Y_{s-},Z_{s-},u)}{Y_{s-}}\right)\nu(du)\bigg\}ds\nonumber\\
&\ge&\lambda,
\end{eqnarray}
then
\begin{eqnarray}\label{per2a}
\liminf_{t\to\infty}\frac{1}{t}\int_0^tY_sds \ge \frac{\lambda}{\lambda_0}\ a.s..
\end{eqnarray}

\noindent (iii) If there exist positive constants $\lambda_0$ and $\lambda$ such that
\begin{eqnarray}\label{per11a}
&&\liminf_{t\to\infty}\inf_{(x,y,z)\in\mathbb{R}^3_+}\bigg\{\lambda_0y+\frac{b_2(t,x,y,z)}{y} - \frac{\sum_{j=1}^n\sigma^2_{2j}(t,x,y,z)}{2y^2}\nonumber\\
	&&\ \ \ \ \ \ \ \ \ \ \ \ \ \ +\int_{\{|u|<1\}}\left[\ln\left(1+\frac{h_2(t,x,y,z,u)}{y}\right)-\frac{h_2(t,x,y,z,u)}{y}\right]\nu(du)\nonumber\\
&&\ \ \ \ \ \ \ \ \ \ \ \ \ \ +\int_{\{|u|\ge1\}}\ln\left(1+\frac{g_2(t,x,y,z,u)}{y}\right)\nu(du)\bigg\}\nonumber\\
&\ge&\lambda,
\end{eqnarray}
then
\begin{eqnarray*}
\liminf_{t\to\infty}\frac{1}{t}\int_0^tY_sds \ge \frac{\lambda}{\lambda_0}\ a.s..
\end{eqnarray*}
\end{thm}

\noindent {\bf Proof.}\ \ (i) By It\^{o}'s formula, we get
\begin{eqnarray}\label{extpf1a}
\ln Y_t&=&\ln Y_0  + \int_0^t \left[ \frac{b_2(s,X_s,Y_s,Z_s)}{Y_s} - \frac{\sum_{j=1}^n\sigma_{2j}^2(s,X_s,Y_s,Z_s)}{2Y_s^2}\right]ds\nonumber\\
&&+\int_0^t\int_{\{|u|<1\}} \left[\ln\left(1+\frac{h_2(s,X_{s-},Y_{s-},Z_{s-},u)}{Y_{s-}}\right)-\frac{h_2(s,X_{s-},Y_{s-},Z_{s-},u)}{Y_{s-}}\right]\nu(du)ds\nonumber\\
&&+\int_0^t\int_{\{|u|\ge1\}} \ln\left(1+\frac{g_2(s,X_{s-},Y_{s-},Z_{s-},u)}{Y_{s-}}\right)\nu(du)ds\nonumber\\
&&+\int_0^t\frac{\sum_{j=1}^n\sigma_{2j}(s,X_s,Y_s,Z_s)}{Y_s}dB_s^{(j)}\nonumber\\
&&+\int_0^t\int_{\{|u|<1\}}\ln\left(1+\frac{h_2(s,X_{s-},Y_{s-},Z_{s-},u)}{Y_{s-}}\right)\tilde{N}(ds,du)\nonumber\\
&&+\int_0^t\int_{\{|u|\ge1\}}\ln\left(1+\frac{g_2(s,X_{s-},Y_{s-},Z_{s-},u)}{Y_{s-}}\right)\tilde{N}(ds,du).
\end{eqnarray}
Denote the martingale part of $\ln Y_t$ by $M_t$. Then, by (\ref{extpf1a}), we get
\begin{eqnarray*}
\langle M\rangle_t&=& \int_0^t \frac{\sum_{j=1}^n\sigma_{2j}^2(s,X_s,Y_s,Z_s)}{Y_s^2}ds\\
&&+\int_0^t\int_{\{|u|<1\}}\left[\ln\left(1+\frac{h_2(s,X_{s-},Y_{s-},Z_{s-},u)}{Y_{s-}}\right)\right]^2\nu(du)ds\nonumber\\
&&+\int_0^t\int_{\{|u|\ge1\}}\left[\ln\left(1+\frac{g_2(s,X_{s-},Y_{s-},Z_{s-},u)}{Y_{s-}}\right)\right]^2\nu(du)ds.
\end{eqnarray*}

By (\ref{p1_eqa}) and the strong law of large numbers for martingales (see \cite[Theorem 10, Chapter 2]{LS}), we get
\begin{eqnarray}\label{extmart1}
\lim\limits_{t\rightarrow\infty}\frac{M_t}{t}=0\ a.s..
\end{eqnarray}
Then,   (\ref{p1_eq9a}) holds by (\ref{extqg}),  (\ref{extpf1a}) and (\ref{extmart1}).

\noindent (ii) By (\ref{per1a}) and (\ref{extpf1a}), if we take $\eta \in (0,\lambda)$ then there exists $T_{\eta} > 0$ such that for $t \ge T_{\eta}$,
\begin{eqnarray*}
\ln Y_t &\ge&\ln Y_0 +  (\lambda-\eta)t -\lambda_0\int_0^tY_sds +\int_0^t\frac{\sum_{j=1}^n\sigma_{2j}(s,X_s,Y_s,Z_s)}{Y_s}dB_t^{(j)}\\
&&+\int_0^t\frac{\sum_{j=1}^n\sigma_{2j}(s,X_s,Y_s,Z_s)}{Y_s}dB_s^{(j)}\nonumber\\
&&+\int_0^t\int_{\{|u|<1\}}\ln\left(1+\frac{h_2(s,X_{s-},Y_{s-},Z_{s-},u)}{Y_{s-}}\right)\tilde{N}(ds,du)\nonumber\\
&&+\int_0^t\int_{\{|u|\ge1\}}\ln\left(1+\frac{g_2(s,X_{s-},Y_{s-},Z_{s-},u)}{Y_{s-}}\right)\tilde{N}(ds,du).
\end{eqnarray*}
Thus, by following the argument of the proof of \cite[Lemma 5.1]{JJ}, we can show that (\ref{per2a}) holds by  (\ref{extmart1}).

\noindent (iii) Obviously, condition (\ref{per11a}) implies condition (\ref{per1a}). Hence, the assertion is a direct consequence of assertion (ii).\hfill\fbox

\begin{rem}
If we take the following assumption of our model:
$$b_2(t,x,y,z) = b_{2,1}(t,x,y,z) - b_{2,2}(t,x,y,z),
$$
where $b_{2,i}(t,x,y,z)\ge 0$ for any $(t,x,y,z)\in[0,\infty)\times \mathbb{R}^3_+$, $i=1,2$, then condition (\ref{extqg}) can be strengthened to
\begin{eqnarray}\label{ext22}
\alpha^*&:=&\limsup\limits_{t\to\infty}\left\{\sup\limits_{(x,y,z)\in\mathbb{R}^3_+}\left[\frac{b^2_{2,1}(t,x,y,z)}{2\sum_{j=1}^n\sigma^2_{2j}(t,x,y,z)}-\frac{b_{2,2}(t,x,y,z)}{y}\right]\right.\nonumber\\
&&\ \ \ \ \ \ \ \ \ \ \ \
+\int_{\{|u|<1\}}\sup\limits_{(x,y,z)\in\mathbb{R}^3_+}\left[\ln\left(1+\frac{h_2(t,x,y,z,u)}{y}\right)-\frac{h_2(t,x,y,z,u)}{y}\right]\nu(du)\nonumber\\
&&\ \ \ \ \ \ \ \ \ \ \ \
+\int_{\{|u|\ge1\}}\sup\limits_{(x,y,z)\in\mathbb{R}^3_+}\left[\ln\left(1+\frac{g_2(t,x,y,z,u)}{y}\right)\right]\nu(du)\nonumber\\
&<& 0,
\end{eqnarray}
In fact, we have $\alpha\le \alpha^*$ and hence condition (\ref{ext22}) implies that
$$
\limsup_{t\rightarrow\infty}\frac{\ln Y_t}{t}\le\alpha^*\ \ a.s..
$$
\end{rem}

Denote by $L^{\infty}_{+}[0,\infty)$ the set of all bounded, non-negative, measurable functions on $[0,\infty)$.  For $f\in L^{\infty}_{+}[0,\infty)$, define
$$
\overline{f}:=\sup_{t\in[0,\infty)}f(t),\ \ \ \ \underline{f}:=\inf_{t\in[0,\infty)}f(t).
$$

\color{black}
\begin{exa}\label{exa3.3}
Let $\Lambda,\mu,\beta,\gamma,\varepsilon,\sigma\in L^{\infty}_{+}[0,\infty)$.
We consider the system
\begin{equation}\label{XC}
\begin{cases} dX_t =[\Lambda(t)-\mu(t) X_t -\beta(t) X_tY_t]dt - \sigma(t) X_tY_tdB_t, \\dY_t = [\beta(t) X_tY_t-(\mu(t)+\gamma(t)+\varepsilon(t))Y_t]dt + \sigma(t) X_tY_tdB_t, \\ dZ_t= [\gamma(t) Y_t-\mu(t) Z_t]dt. \end{cases}
\end{equation}

Suppose that
$$
\underline{\mu}>0.
$$
By (\ref{XC}), we get
$$
d(X_t+Y_t+Z_t)\le [\overline{\Lambda}-\underline{\mu}(X_t+Y_t+Z_t)]dt,
$$
which implies that
$$
\Gamma:=\left\{(x,y,z)\in \mathbb{R}^3_+:x+y+z\le \frac{\overline{\Lambda}}{\underline{\mu}}\right\}
$$
is an invariant set of the system (\ref{XC}). Hence,  the system (\ref{XC}) has a unique strong solution taking values in $\Gamma$ by Theorem \ref{RE}.

Define
$$
\alpha^*:=\sup_{x\in\left(0,\frac{\overline{\Lambda}}{\underline{\mu}}\right)}\left[\overline{\beta} x-(\underline{\mu}+\underline{\gamma}+\underline{\varepsilon})-\frac{\underline{\sigma}^2 x^2}{2}\right].
$$
Then, we have that
\begin{eqnarray*}
&&\text{ condition (\ref{extqg})} \\
&\Leftrightarrow& \alpha = \limsup_{t\rightarrow\infty}\sup_{x\in\left(0,\frac{\overline{\Lambda}}{\underline{\mu}}\right)}\left[\beta(t) x-(\mu(t)+\gamma(t)+\varepsilon(t))-\frac{\sigma^2(t) x^2}{2}\right]<0\\
&\Leftarrow&\alpha \le\alpha^*<0\\
&\Leftrightarrow&\alpha^*=\max\left\{\frac{\overline{\beta}\,\overline{\Lambda}}{\underline{\mu}}-(\underline{\mu}+\underline{\gamma}
+\underline{\varepsilon})-\frac{\underline{\sigma}^2\overline{\Lambda}^2}{2\underline{\mu}^2},\ \frac{\overline{\beta}^2}{2\underline{\sigma}^2}
-(\underline{\mu}+\underline{\gamma}+\underline{\varepsilon})\right\}<0\\
&\Leftrightarrow&\left\{\begin{array}{ll}
\alpha^*=\frac{\overline{\beta}\,\overline{\Lambda}}{\underline{\mu}}-(\underline{\mu}+\underline{\gamma}+\underline{\varepsilon})
-\frac{\underline{\sigma}^2\overline{\Lambda}^2}{2\underline{\mu}^2}<0, &{\rm if}\ \underline{\sigma}^2\le\frac{\underline{\mu}\overline{\beta}}{\overline{\Lambda}},\\
\alpha^*=\frac{\overline{\beta}^2}{2\underline{\sigma}^2}
-(\underline{\mu}+\underline{\gamma}+\underline{\varepsilon})<0, &{\rm if}\ \underline{\sigma}^2>\frac{\underline{\mu}\overline{\beta}}{\overline{\Lambda}}. \end{array}\right.
\end{eqnarray*}
Thus, by Theorem \ref{p1_thm2}(i), we obtain that if
$$
\underline{\sigma}^2\le\frac{\underline{\mu}\overline{\beta}}{\overline{\Lambda}}\ \ {\rm and}\ \ \tilde{R}_0:=\frac{\overline{\beta}\,\overline{\Lambda}}{\underline{\mu}(\underline{\mu}+\underline{\gamma}+
\underline{\varepsilon})}-\frac{\underline{\sigma}^2\overline{\Lambda}^2}{2\underline{\mu}^2(\underline{\mu}+\underline{\gamma}+\underline{\varepsilon})}<1,
$$
then the disease extincts with exponential rate
$$
-\alpha\ge (\underline{\mu}+\underline{\gamma}+\underline{\varepsilon})\left(1-\tilde{R}_0\right);
$$
if
$$
\underline{\sigma}^2>\max\left\{\frac{\underline{\mu}\overline{\beta}}{\overline{\Lambda}},\ \frac{\overline{\beta}^2}{2(\underline{\mu}
+\underline{\gamma}+\underline{\varepsilon})}\right\},
$$
then the disease extincts with exponential rate
$$
-\alpha\ge(\underline{\mu}+\underline{\gamma}+\underline{\varepsilon})-\frac{\overline{\beta}^2}{2\underline{\sigma}^2}.
$$
This result generalizes the result given in \cite[Theorem 2.1]{JJ}.

By (\ref{XC}), we get
$$
X_t+Y_t=X_0+Y_0+\int_0^t[\Lambda(s)-\mu(s) X_s-(\mu(s)+\gamma(s)+\varepsilon(s))Y_s]ds.
$$
Since
\begin{equation}\label{RRR}
X_t+Y_t\le \frac{\overline{\Lambda}}{\underline{\mu}},\ \ \ \ \forall t\ge 0,
\end{equation}
we get
$$
\lim_{t\rightarrow\infty}\frac{1}{t}\int_0^t[\Lambda(s)-\mu(s) X_s-(\mu(s)+\gamma(s)+\varepsilon(s))Y_s]ds=0,
$$
which implies that
\begin{equation}\label{R1}
\lim_{t\rightarrow\infty}\frac{1}{t}\int_0^tX_sds\ge\frac{\underline{\Lambda}}{\overline{\mu}}-\lim_{t\rightarrow\infty}\frac{\overline{\mu}+\overline{\gamma}
+\overline{\varepsilon}}{\overline{\mu} t}\int_0^tY_sds.
\end{equation}

Suppose that
$$
\tilde{R}_0:=\frac{\underline{\beta}\underline{\Lambda}}{\overline{\mu}(\overline{\mu}+\overline{\gamma}+\overline{\varepsilon})}
-\frac{\overline{\sigma}^2\overline{\Lambda}^2}{2\underline{\mu}^2(\overline{\mu}+\overline{\gamma}+\overline{\varepsilon})}>1.
$$
Then, by (\ref{XC})--(\ref{R1}), we get
\begin{eqnarray*}
&&\liminf_{t\rightarrow\infty}\frac{1}{t}\int_0^t\left[\frac{\underline{\beta}(\overline{\mu}+\overline{\gamma}
+\overline{\varepsilon})}{\overline{\mu}}\cdot Y_s+\frac{b_2(s,X_s,Y_s,Z_s)}{Y_s}
-\frac{\sigma^2(s)X_s^2}{2}\right]ds\\
&\ge&\frac{\underline{\beta}\underline{\Lambda}}{\overline{\mu}}-(\overline{\mu}+\overline{\gamma}+\overline{\varepsilon})-\frac{\overline{\sigma}^2\overline{\Lambda}^2}{2\underline{\mu}^2}.
\end{eqnarray*}
Therefore, by Theorem \ref{p1_thm2}(ii), we obtain that the disease is persistent and
$$
\liminf_{t\rightarrow\infty}\frac{1}{t}\int_0^tY_sds\ge \frac{\overline{\mu}(\tilde{R}_0-1)}{\underline{\beta}}\ \ a.s..
$$
This result generalizes the result given in  \cite[Theorem 3.1]{JJ}.
\end{exa}

Let $M>0$ be a fixed constant. For $x\ge 0$, define
$$
x^{\star} := x \wedge 1,\ \ \ \  x^\dagger  := x \wedge M.
$$
\begin{exa}\label{exa3.4}
In the following examples, we let $d=1$ and the intensity measure $\nu$ of the Poisson random measure $N$ be given by
$$
d\nu =\mathds{1}_{[-2,2]}(x)dx,
$$
where $dx$ is the Lebesgue measure.

\noindent (a) Let $\Lambda, \mu, \beta, \gamma_1,\gamma_2,\gamma_3,\gamma_4,\xi, \sigma_1, \sigma_2 , \varphi_1,\varphi_2,\varphi_3 \in  L^{\infty}_{+}[0,\infty)$ and $h_1, h_2, h_3, g_1,g_2 \in L^\infty_+(-\infty,\infty)$. Define
$$
\varphi(t,x,y) = \varphi_1(t)x+ \varphi_2(t)y+ \varphi_3(t)xy,\ \ \ \ (t,x,y)\in[0,\infty)\times \mathbb{R}^2_+.
$$
We consider the system
\begin{eqnarray}\label{ex34a}
\begin{cases}
\vspace{.35cm}
dX_t =\left[\Lambda(t) - \mu(t) \xd  -\frac{\beta(t)\xd{\yd}^{\xi(t)}}{1+\varphi(t,X_t,Y_t)} + \gamma_1(t)\zd\right]dt - \frac{\sigma_1(t)\xd{\yd}^{\xi(t)}}{1+\varphi(t,X_t,Y_t)} dB_t^{(1)}\\ \vspace{.35cm}
\indent\indent- \int_{\{|u| < 1\}} [h_1(u)\xms\yms-h_3(u)\xms\zms ]\tilde{N}(dt,du)
- \int_{\{|u| \ge 1\}} g_1(u)\xms\yms N(dt,du), \\

\vspace{.35cm}
dY_t = \left [ \frac{\beta(t)\xd{\yd}^{\xi(t)}}{1+\varphi(t,X_t,Y_t)} + (\gamma_2(t)-\mu(t)-\gamma_3(t)\yd)\yd \right]dt +   \frac{\sigma_1(t)\xd{\yd}^{\xi(t)}}{1+\varphi(t,X_t,Y_t)} dB_t^{(1)} + \sigma_2(t)\yd\zd dB_t^{(2)} \\ \vspace{.35cm}
\indent\indent + \int_{\{|u| < 1\}} [h_1(u)\xms\yms -h_2(u)\yms\zms]\tilde{N}(dt,du) \\ \vspace{.35cm}
\indent\indent +\int_{\{|u| \ge 1\}}[ g_1(u)\xms\yms -g_2(u)\yms\zms ]N(dt,du), \\
\vspace{.35cm}
dZ_t = \left[\gamma_4(t)\yd  - (\mu(t)+\gamma_1(t))\zd\right]dt - \sigma_2(t)\yd\zd dB_t^{(2)}
	\\ \vspace{.35cm} \indent\indent+ \int_{\{|u| < 1\}}[ h_2(u)\yms\zms -h_3(u)\xms\zms] \tilde{N}(dt,du) +\int_{\{|u| \ge 1\}}g_2(u)\yms\zms N(dt,du).
\end{cases}
\end{eqnarray}

Suppose that
$$
\underline{\xi}\ge1,\
\overline{h_i},\overline{g_j}<1,\ \ \ \ i=1,2,3,\  j=1,2.
$$
Then, Assumptions \ref{assm21}--\ref{assm24} hold. Thus, by Theorem \ref{RE}, the system  (\ref{ex34a}) has a unique strong solution in $\mathbb{R}^3_+$.
Assume that
$$
\overline{\mu}<\underline{\gamma_2},\ \ \ \ (M\vee 1)^{2\overline{\xi}}(\overline{\sigma_1}^2 + \overline{\sigma_2}^2) +4[\overline{h_1} -  \ln\{(1- \overline{h_2})(1- \overline{g_2}) \}] < 2\min\{M,\underline{\gamma_2}-\overline{\mu}\}.
$$
Set
$$
\lambda_0 =\overline{\gamma_3}+1,\ \ \ \ \lambda=\min\{M,\underline{\gamma_2}-\overline{\mu}\}-\left\{\frac{(M\vee 1)^{2\overline{\xi}}(\overline{\sigma_1}^2 + \overline{\sigma_2}^2)}{2} +2[\overline{h_1} -  \ln\{(1- \overline{h_2})(1- \overline{g_2}) \}]\right\}.
$$
Then, by Theorem \ref{p1_thm2}(iii), we obtain that the disease is persistent and
$$
	\liminf_{t \to \infty} \frac{1}{t}\int_0^tY_sds \ge\frac{\lambda}{\lambda_0} \ \  a.s..
$$

\noindent (b) Let $\Lambda, \mu, \beta, \gamma_1,\gamma_2, \sigma \in L^{\infty}_+[0,\infty)$ and $h_1, h_2, h_3, g_1,g_2 \in L^\infty_+(-\infty,\infty)$. We consider the system
\begin{eqnarray}\label{ex34b}
\begin{cases}
\vspace{.35cm}
dX_t =\left[\Lambda(t) -\mu(t)\xd -\beta(t) \xd\yd+\gamma_1(t)\zd\right]dt - \sigma(t)\xd\yd\zd dB_t \\ \vspace{.35cm}\indent\indent
- \int_{\{|u| < 1\}} [h_1(u)-h_3(u)]\xms\yms\zms\tilde{N}(dt,du)\\ \vspace{.35cm}\indent\indent
- \int_{\{|u| \ge 1\}}[ g_1(u)-g_3(u)]\xms\yms\zms N(dt,du), \\

\vspace{.35cm}
dY_t = \left[\beta(t)\xd\yd-(\mu(t)+\gamma_2(t))\yd\right]dt+2\sigma(t) \xd\yd\zd dB_t  \\ \vspace{.35cm}
\indent\indent + \int_{\{|u| < 1\}} [h_1(u)-h_2(u)]\xms\yms\zms\tilde{N}(dt,du) \\ \vspace{.35cm}
\indent\indent +\int_{\{|u| \ge 1\}}[ g_1(u)-g_2(u)]\xms\yms\zms N(dt,du), \\
\vspace{.35cm}
dZ_t= [\gamma_2(t)\yd-(\mu(t)+\gamma_1(t))\zd ]dt-  \sigma(t) \xd\yd\zd dB_t
\\ \vspace{.35cm}\indent\indent+ \int_{\{|u| < 1\}} [h_2(u)-h_3(u)]\xms\yms\zms\tilde{N}(dt,du) \\ \vspace{.35cm}
 \indent\indent+\int_{\{|u| \ge 1\}}[g_2(u)-g_3(u)]\xms\yms\zms N(dt,du).
\end{cases}
\end{eqnarray}
\end{exa}

Suppose that
$$
\overline{h_i},\overline{g_j}<1,\ \ \ \ i=1,2,3,\ j=1,2.
$$
Then, Assumptions \ref{assm21}--\ref{assm24} hold. Thus, by Theorem \ref{RE}, the system (\ref{ex34b}) has a unique strong solution in $\mathbb{R}_+^3$. If
$$
 \overline{\beta} +2\overline{g_1} < \underline{\gamma_2} + \underline{\mu},
$$
then by Theorem \ref{p1_thm2}(i) and noting that $\ln(1+x)-x\le 0$ for $x>-1$, we obtain that
the disease extincts with exponential rate
$$
 - \alpha >  \underline{\gamma_2} + \underline{\mu} - \overline{\beta} - 2\overline{g_1}.
$$

\section{Model for population proportions}\setcounter{equation}{0}

In this section, we let $X_t$, $Y_t$ and $Z_t$ denote respectively the proportions of susceptible, infected and recovered populations at time $t$. This is a special case of the USSIR  model (\ref{p1_eq2}) which has been considered in the literature (cf. \cite{EL}, \cite{G} and \cite{T}). Define
$$
\Delta:=\{(x,y,z)\in \mathbb{R}^3_+:x+y+z=1\}.
$$

We make the following assumptions.
\begin{enumerate}[\bf({B}1)]

\item\label{assm1} There exists $(x_0,y_0,z_0)\in\Delta$ such that for any $T \in (0,\infty)$ and $i=1,2,3$,
$$
	b_i(\cdot,x_0,y_0,z_0),\ \ \sum\limits_{j=1}^n|\sigma_{ij}(\cdot,x_0,y_0,z_0)| \in L^2[0,T],\ \ \int_{\{|u|<1\}} |h_i(\cdot,x_0,y_0,z_0,u)|^2\nu(du) \in L^1[0,T].
$$
\item\label{assm2} For any $T \in (0,\infty)$ and $N \in \mathbb{N}$, there exists $K_{N,T} \in L_+^1[0,T]$ such that
\begin{eqnarray*}
&&\sum\limits_{i=1}^3  \lvert b_i(t,x_1,y_1,z_1)-b_i(t,x_2,y_2,z_2)\rvert^2 +\sum\limits_{i=1}^3\sum\limits_{j=1}^n\lvert \sigma_{ij}(t,x_1,y_1,z_1)-\sigma_{ij}(t,x_2,y_2,z_2)\rvert^2\\
	&&+\sum\limits_{i=1}^3\int_{\{|u|<1\}}| h_i(t,x_1,y_1,z_1,u) - h_i(t,x_2,y_2,z_2,u)|^2\nu(du)\\
& \le& K_{N,T}(t)[(x_1-x_2)^2+(y_1-y_2)^2+(z_1-z_2)^2],\\
&&\forall t\in [0,T],\ (x_1,y_1,z_1),(x_2,y_2,z_2)\in \left[\frac{1}{N},1-\frac{1}{N}\right]^3.
\end{eqnarray*}

\item\label{assm3}For any $t\in(0,\infty)$, $(x,y,z) \in \Delta$  and $u\in \mathbb{R}^d-\{0\}$,
\begin{eqnarray*}
&&\sum\limits_{i=1}^3 b_i(t,x,y,z)=0, \ \ \sum\limits_{i=1}^3 \sigma_{ij}(t,x,y,z)=0\text{ for } j=1,2,\ldots, n,\\
	&&  \sum\limits_{i=1}^3h_i(t,x,y,z,u) = 0,\ \   \sum\limits_{i=1}^3g_i(t,x,y,z,u) = 0.
\end{eqnarray*}

\item\label{assm4} For any $(x,y,z) \in \Delta$, $t\in(0,\infty)$ and $u\in \mathbb{R}^d-\{0\}$,
\begin{eqnarray*}
&&\left( 1+ \frac{h_1(t,x,y,z,u)}{x}\right),\ \left( 1+ \frac{h_2(t,x,y,z,u)}{y}\right),\ \left( 1+ \frac{h_3(t,x,y,z,u)}{z}\right) > 0,
\end{eqnarray*}
and
\begin{eqnarray*}
&&\left( 1+ \frac{g_1(t,x,y,z,u)}{x}\right),\ \left( 1+ \frac{g_2(t,x,y,z,u)}{y}\right),\ \left( 1+ \frac{g_3(t,x,y,z,u)}{z}\right) > 0.
\end{eqnarray*}

\item\label{assm5} For any $T\in (0,\infty)$,
\begin{eqnarray*}
&&\inf_{t\in[0, T],\, (x,y,z)\in\Delta}\left\{\frac{b_1(t,x,y,z)}{x}+\frac{b_2(t,x,y,z)}{y}+\frac{b_3(t,x,y,z)}{z}\right\}>-\infty,
\end{eqnarray*}
\begin{eqnarray*}
&&\sum_{j=1}^n\sup_{t\in[0, T],\, (x,y,z)\in\Delta}\left\{\frac{|\sigma_{1j}(t,x,y,z)|}{x}+\frac{|\sigma_{2j}(t,x,y,z)|}{y}+\frac{|\sigma_{3j}(t,x,y,z)|}{z}\right\}<\infty,
\end{eqnarray*}
and there exists $\eta_T\in {\cal B}(\mathbb{R}^d)$ such that
(\ref{p1_eq4}) holds and
$$
\sup_{t\in[0, T],\, (x,y,z)\in\Delta, u\in \mathbb{R}^d-\{0\}}\frac{k(t,x,y,z,u)}{\eta_T(u)}<\infty,
$$
where $k(t,x,y,z,u)$ is defined by (\ref{p1_eq3}).
\end{enumerate}

We now discuss the existence and uniqueness of solutions to the system (\ref{p1_eq2}) when considering the proportional form.
\begin{thm}\label{thm0}
Suppose that Assumptions \ref{assm1}--\ref{assm5} hold. Then, for any given initial value $(X_0,Y_0, Z_0) \in \Delta$, the system (\ref{p1_eq2}) has a unique strong solution taking values in $\Delta$.
\end{thm}

\noindent {\bf Proof.}\ \
By \ref{assm3}, \ref{assm4} and the interlacing technique, to complete the proof, we need only  consider the case that $g_i\equiv0$, $i=1,2,3$. Then, equation (\ref{p1_eq2}) becomes  equation (\ref{p1_eq5}).

By \ref{assm1} and \ref{assm2}, similar to \cite[Lemma 2.1]{Sun}, we can show that there exists a unique local strong solution to equation (\ref{p1_eq5}) on $[0,\tau)$, where $\tau$ is the explosion time. We will show below that $\tau = \infty$ a.s.. Define
\begin{eqnarray*}
\tau_N =\inf\left\{t\in[0,\tau): (X_t,Y_t,Z_t)\not\in\left[\frac{1}{N},1-\frac{1}{N}\right]^3\right\},\ \ \ \ N\in\mathbb{N},
\end{eqnarray*}
and
$$
\tau_{\infty}=\lim_{N\rightarrow\infty}\tau_N.
$$
We have that $\tau_{\infty} \le \tau$ so it suffices to show $\tau_{\infty} = \infty$ a.s.. Hence assume the contrary that there exist $\varepsilon>0$ and $T>0$ such that
$$
\mathbb{P}(\tau_{\infty} < T) >\varepsilon,
$$
which implies that
\begin{eqnarray}\label{p1_eq6}
\mathbb{P}(\tau_N < T) >\varepsilon,\ \ \ \ \forall N\in\mathbb{N}.
\end{eqnarray}

Define
$$
V(x,y,z) = -\ln(xyz),\ \ \ \ (x,y,z) \in (0,1)^3,
$$
and
$$
W_t=(X_t,Y_t,Z_t).
$$
By It\^{o}'s formula, we obtain that for $t\le \tau_N$,
\begin{eqnarray*}
V(W_t)&=& V(W_0)+\int_0^t\mathcal{L}V(s,W_s)ds +\int_0^t \langle V_x(W_s), \sigma(s,X_s) \rangle dB_s\\
&&+ \int_0^t\int_{\{|u|<1\}}\left[V(W_{s-}+h(s,W_{s-},u)) - V(W_{s-}) \right]\tilde{N}(ds,du).
\end{eqnarray*}
 Then, by \ref{assm5}, there exist $C_T>0$ and  $\eta_T\in {\cal B}(\mathbb{R}^d)$ such that (\ref{p1_eq4}) holds and
\begin{eqnarray*}
&&(\ln N)\mathbb{P}(\tau_N < T) -V(X_0,Y_0,Z_0)\\
&\le&\mathbb{E}\left[V(X_{T\wedge \tau_N},Y_{T\wedge \tau_N},Z_{T\wedge \tau_N})\right] -V(X_0,Y_0,Z_0)\\
&=&\mathbb{E}\left[\int_0^{T\wedge \tau_N} \mathcal{L}V(s,X_s,Y_s,Z_s) ds\right]\\
&=&- \mathbb{E}\left[\int_0^{T\wedge \tau_N} \left\{\frac{b_1(s,X_s,Y_s,Z_s)}{X_s}+\frac{b_2(s,X_s,Y_s,Z_s)}{Y_s}+\frac{b_3(s,X_s,Y_s,Z_s)}{Z_s} \right\}ds\right]\\
	&&+\ \frac{1}{2}\sum\limits_{j=1}^n\mathbb{E}\left[\int_0^{T\wedge \tau_N} \left\{\frac{\sigma_{1j}^2(s,X_s,Y_s,Z_s)}{X_s^2}+\frac{\sigma_{2j}^2(s,X_s,Y_s,Z_s)}{Y_s^2}+\frac{\sigma_{3j}^2(s,X_s,Y_s,Z_s)}{Z_s^2}\right\}ds\right]\\
	&&+\mathbb{E}\left[\int_0^{T\wedge \tau_N}\ \int_{\{|u|<1\}}k(s,X_s,Y_s,Z_s,u)\nu(du)ds\right]\\
&\le& C_TT+\frac{C_T^2T}{2}+C_TT \int_{\{|u|<1\}}\eta_T(u)\nu(du),
\end{eqnarray*}
which contradicts with (\ref{p1_eq6}). Therefore, $\tau = \infty$ a.s. and the proof is complete. \hfill\fbox


Similar to Theorem \ref{p1_thm2}, we can prove the following result on the extinction and persistence of diseases.

\begin{thm}\label{p1_thm1}
Suppose that Assumptions \ref{assm1}--\ref{assm5} hold. Let $(X_t,Y_t,Z_t)$ be a solution to equation (\ref{p1_eq2}) with $(X_0,Y_0,Z_0) \in \Delta$. We assume that
$$
\int_0^{\infty}\frac{\varphi(t)}{(1+t)^2}dt<\infty,
$$
where
\begin{eqnarray*}
\varphi(t)&:=&\sup_{(x,y,z) \in \Delta} \left\{\frac{ \sum_{j=1}^n \sigma^2_{2j}(t,x,y,z)}{y^2}+\int_{\{|u|<1\}}\left[\ln\left(1+\frac{h_2(t,x,y,z,u)}{y}\right)\right]^2\nu(du)\right.\\
&&\ \ \ \ \ \ \ \ \ \ \ \left.+\int_{\{|u|\ge1\}}\left[\ln\left(1+\frac{g_2(t,x,y,z,u)}{y}\right)\right]^2\nu(du)\right\}.
\end{eqnarray*}
\noindent (i) If
\begin{eqnarray*}
\alpha&:=&\limsup\limits_{t\to\infty}\left\{\sup\limits_{(x,y,z)\in\Delta}\left[\frac{b_2(t,x,y,z)}{y} - \frac{\sum_{j=1}^n\sigma^2_{2j}(t,x,y,z)}{2y^2}\right]\right.\nonumber\\
&&\ \ \ \ \ \ \ \ \ \ \ \
+\int_{\{|u|<1\}}\sup\limits_{(x,y,z)\in\Delta}\left[\ln\left(1+\frac{h_2(t,x,y,z,u)}{y}\right)-\frac{h_2(t,x,y,z,u)}{y}\right]\nu(du)\nonumber\\
&&\ \ \ \ \ \ \ \ \ \ \ \
+\int_{\{|u|\ge1\}}\sup\limits_{(x,y,z)\in\Delta}\left[\ln\left(1+\frac{g_2(t,x,y,z,u)}{y}\right)\right]\nu(du)\nonumber\\
&<& 0,
\end{eqnarray*}
then
\begin{eqnarray*}
\limsup_{t\to\infty} \frac{\ln Y_t}{t} \le \alpha\ a.s..
\end{eqnarray*}

\noindent (ii) If there exist positive constants $\lambda_0$ and $\lambda$ such that
\begin{eqnarray*}
&&\liminf_{t\to\infty}\frac{1}{t}\int_0^t\bigg\{\lambda_0Y_s+\frac{b_2(s,X_s,Y_s,Z_s)}{Y_s} - \frac{\sum_{j=1}^n\sigma^2_{2j}(s,X_s,Y_s,Z_s)}{2Y_s^2}\nonumber\\
	&&\ \ \ \ \ \ \ \ \ \ \ \ \ \ +\int_{\{|u|<1\}}\left[\ln\left(1+\frac{h_2(s,X_{s-},Y_{s-},Z_{s-},u)}{Y_{s-}}\right)-\frac{h_2(s,X_{s-},Y_{s-},Z_{s-},u)}{Y_{s-}}\right]\nu(du)\nonumber\\
&&\ \ \ \ \ \ \ \ \ \ \ \ \ \ +\int_{\{|u|\ge1\}}\ln\left(1+\frac{g_2(s,X_{s-},Y_{s-},Z_{s-},u)}{Y_{s-}}\right)\nu(du)\bigg\}ds\nonumber\\
&\ge&\lambda,
\end{eqnarray*}
then
\begin{eqnarray*}
\liminf_{t\to\infty}\frac{1}{t}\int_0^tY_sds \ge \frac{\lambda}{\lambda_0}\ a.s..
\end{eqnarray*}

\noindent (iii) If there exist positive constants $\lambda_0$ and $\lambda$ such that
\begin{eqnarray}\label{per1188}
&&\liminf_{t\to\infty}\inf_{(x,y,z)\in\Delta}\bigg\{\lambda_0y+\frac{b_2(t,x,y,z)}{y} - \frac{\sum_{j=1}^n\sigma^2_{2j}(t,x,y,z)}{2y^2}\nonumber\\
	&&\ \ \ \ \ \ \ \ \ \ \ \ \ \ +\int_{\{|u|<1\}}\left[\ln\left(1+\frac{h_2(t,x,y,z,u)}{y}\right)-\frac{h_2(t,x,y,z,u)}{y}\right]\nu(du)\nonumber\\
&&\ \ \ \ \ \ \ \ \ \ \ \ \ \ +\int_{\{|u|\ge1\}}\ln\left(1+\frac{g_2(t,x,y,z,u)}{y}\right)\nu(du)\bigg\}\nonumber\\
&\ge&\lambda,
\end{eqnarray}
then
\begin{eqnarray*}
\liminf_{t\to\infty}\frac{1}{t}\int_0^tY_sds \ge \frac{\lambda}{\lambda_0}\ a.s..
\end{eqnarray*}
\end{thm}

\color{black}
\begin{exa}\label{exa2.5} We revisit Example \ref{exa3.4} with some changes for the population proportions model. Let $d=1$ and the intensity measure $\nu$ of the Poisson random measure $N$ be given by
$$
d\nu =\mathds{1}_{[-2,2]}(x)dx.
$$

\noindent (a) Let $ \beta, \gamma,\xi, \sigma_1, \sigma_2 , \varphi_1,\varphi_2,\varphi_3 \in  L^{\infty}_{+}[0,\infty)$ and $h_1,h_2, g_1,g_2 \in L^\infty_+(-\infty,\infty)$. Define
$$
\varphi(t,x,y) = \varphi_1(t)x+ \varphi_2(t)y+ \varphi_3(t)xy,\ \ \ \ (t,x,y)\in[0,\infty)\times \mathbb{R}^2_+.
$$
We consider the system
\begin{eqnarray}\label{p1_ex1}
\begin{cases}
\vspace{.35cm}
dX_t = -\frac{\beta(t)X_t^{\xi(t)}Y_t}{1+\varphi(t,X_t,Y_t)}dt - \frac{\sigma_1(t)X_t^{\xi(t)}Y_t}{1+\varphi(t,X_t,Y_t)} dB_t^{(1)}
- \int_{\{|u| < 1\}} h_1(u)X_{t-}^{\xi(t)}Y_{t-}\tilde{N}(dt,du)   \\ \vspace{.35cm}
\indent\indent- \int_{\{|u| \ge 1\}} g_1(u)X_{t-}^{\xi(t)}Y_{t-}N(dt,du), \\

\vspace{.35cm}
dY_t = \left [ \frac{\beta(t)X_t^{\xi(t)}Y_t}{1+\varphi(t,X_t,Y_t)} - \gamma(t)Y_t \right]dt +   \frac{\sigma_1(t)X_t^{\xi(t)}Y_t}{1+\varphi(t,X_t,Y_t)} dB_t^{(1)} + \sigma_2(t)Y_tdB_t^{(2)} \\ \vspace{.35cm}
\indent\indent + \int_{\{|u| < 1\}} [h_1(u)X_{t-}^{\xi(t)}Y_{t-}-h_2(u)Y_{t-}]\tilde{N}(dt,du) \\ \vspace{.35cm}
\indent\indent +\int_{\{|u| \ge 1\}}[ g_1(u)X_{t-}^{\xi(t)}Y_{t-}-g_2(u)Y_{t-}]N(dt,du), \\
\vspace{.35cm}
dZ_t = \gamma(t)Y_tdt - \sigma_2(t)Y_tdB_t^{(2)}+ \int_{\{|u| < 1\}} h_2(u)Y_{t-}\tilde{N}(dt,du) \\ \vspace{.35cm} \indent\indent+\int_{\{|u| \ge 1\}}g_2(u)Y_{t-}N(dt,du).
\end{cases}
\end{eqnarray}

Suppose that
$$
\underline{\xi}\ge1,\ \
\overline{h_1},\overline{h_2},\overline{g_1},\overline{g_2}<1.
$$
We have $\frac{dX_t}{dt} + \frac{dY_t}{dt} + \frac{dZ_t}{dt} = 0$. Hence, by Theorem \ref{thm0}, the system (\ref{p1_ex1}) has a unique strong solution taking values in $\Delta$. If
$$
	\overline{\beta}+2\overline{g_1}<\underline{\gamma},
$$
then by Theorem \ref{p1_thm1}(i) and noting that $\ln(1+x)-x\le 0$ for $x>-1$, we obtain that
the disease extincts with exponential rate
$$
	-\alpha \ge \underline{\gamma}-\overline{\beta}-2\overline{g_1}.
$$
Additionally, a key feature of the system (\ref{p1_ex1}) to note is that the transmission function is in the form of power function which differs from the often seen bilinear form.
\color{black}

\newpage
\noindent (b) Let $\beta, \gamma_1,\gamma_2, \sigma_1,\sigma_2,\sigma_3 \in L^{\infty}_+[0,\infty)$ and $h_1,h_2,h_3, g_1,g_2,g_3 \in L^\infty_+(-\infty,\infty)$. We consider the system
\begin{eqnarray}\label{ex1b}
\begin{cases}
\vspace{.35cm}
dX_t = -\beta(t) X_tY_tdt - \sigma_1(t) X_tY_tdB_t
- \int_{\{|u| < 1\}} h_1(u)X_{t-}Y_{t-}\tilde{N}(dt,du)   \\ \vspace{.35cm}
\indent\indent- \int_{\{|u| \ge 1\}} g_1(u)X_{t-}Y_{t-}N(dt,du), \\

\vspace{.35cm}
dY_t = \left[\beta(t) X_t-\gamma_1(t)+\gamma_2(t)Z_t\right]Y_tdt+[\sigma_1(t) X_t-\sigma_2(t) +\sigma_3(t)Z_t]Y_tdB_t \\ \vspace{.35cm}
\indent\indent + \int_{\{|u| < 1\}} [h_1(u)X_{t-}-h_2(u)-h_3(u)Z_{t-}]Y_{t-}\tilde{N}(dt,du) \\ \vspace{.35cm}
\indent\indent +\int_{\{|u| \ge 1\}}[g_1(u)X_{t-}-g_2(u)-g_3(u)Z_{t-}]Y_{t-}N(dt,du), \\
\vspace{.35cm}
dZ_t= [\gamma_1(t)-\gamma_2(t)Z_t]Y_tdt + [\sigma_2(t) - \sigma_3(t)Z_t]Y_tdB_t+\int_{\{|u| < 1\}} [h_2(u)+h_3(u)Z_{t-}]Y_{t-}\tilde{N}(dt,du)\\ \vspace{.35cm} \indent\indent+\int_{\{|u| \ge 1\}}[g_2(u)+g_3(u)Z_{t-}]Y_{t-}N(dt,du).
\end{cases}
\end{eqnarray}
We have $\frac{dX_t}{dt} + \frac{dY_t}{dt} + \frac{dZ_t}{dt} = 0$. Hence, by Theorem \ref{thm0}, the system (\ref{ex1b}) has a unique strong solution taking values in $\Delta$.

Suppose that
$$
\overline{h_1},\,\overline{h_2}+\overline{h_3},\,\overline{g_1},\,\overline{g_2}+\overline{g_3}<1,\ \ \  \overline{\gamma_1}< \underline{\beta} \le \underline{\gamma_2},
$$
and
$$
[(\overline{\sigma_1}+\overline{\sigma_3})\vee\overline{\sigma_2}]^2+4[\overline{h_1}-\ln\{(1-\overline{h_2}-\overline{h_3})(1-\overline{g_2}-\overline{g_3})\}] < 2(\underline{\gamma_2} - \overline{\gamma_1}).
$$
Set
$$\lambda_0=\underline{\gamma_2},\ \ \ \ \lambda=\underline{\gamma_2}-\overline{\gamma_1} - \left\{\frac{[(\overline{\sigma_1}+\overline{\sigma_3})\vee\overline{\sigma_2}]^2}{2}+2[\overline{h_1}-\ln\{(1-\overline{h_2}-\overline{h_3})(1-\overline{g_2}-\overline{g_3})\}]\right\}.
$$ Then, condition (\ref{per1188}) is satisfied. Therefore, by Theorem \ref{p1_thm1}(iii), we obtain that the disease is persistent and
$$
\liminf_{t\rightarrow\infty}\frac{1}{t}\int_0^tY_sds\ge \frac{\lambda}{\lambda_0}\ \ a.s..
$$
\end{exa}

\section{Simulations}\setcounter{equation}{0}

In this section, we present simulations corresponding to Examples  \ref{exa3.3}, \ref{exa3.4} and \ref{exa2.5}. Simulations are completed by use of the {\it Julia} programming language (cf. \cite{julia}) and the {\it DifferentialEquations.jl} package (cf. \cite{DiffEq}). The methodology used is a jump-adapted Euler-Maruyama (EM) scheme (cf. \cite{sdej}) with a time step $\Delta t = 0.01$. We include both average and sample paths in the simulations in the following simulation study. Moreover, the average was computed from $100$ simulations; and whence, $3$ sample paths were selected randomly and plotted with the average. For the readers who seek more information about the {\it Julia} programming language and/or simulation of stochastic differential equations, we refer them to the above references.

\begin{rem}
The time $t$ in the following is taken to be epidemiological time without specific unit; however, we may imagine the time units represent days, weeks or months. Additionally, the evenly spaced-timestep $\Delta t$ corresponds to the hypothetical times at which measurements were taken corresponding to the model; namely, the hypothetical time-series data. Certainly, when analyzing real-world data a specified time measurement would be given.

The choice's made by the authors of parameter values and parameter periods were done so arbitrarily yet educated. Namely, it is important to show the flexibility of the USSIR model and illustrate the viability of the theoretical results by use of simulations. The results are intended to be both understandable and useful to the informed reader.
\end{rem}

\begin{simu} This simulation is concerned with Example \ref{exa3.3}, the system (\ref{XC}). The initial condition is set to $(X_0,Y_0,Z_0)=\left(2.0, 0.8, 1\right)$, where the starting population is $3.8$ million. In the following simulations, the parameters will change to demonstrate their effects on a system with unchanging initial condition. The first two simulations illustrate extinction of the disease and the final simulation will illustrate persistence of the disease. We initially set the parameters in Table 1:

\begin{center}
	\resizebox{9cm}{!}{
    \begin{tabu}{| c  | c | c |}
        \hline
    $f(t)$ &  $\underline{f}$ & $\overline{f}$\\ \hline
    $\beta(t)=0.13+0.01\sin(t)$ & $ 0.12$  & $ 0.14$ \\
    \hline
    $\gamma(t) = 0.9+0.02\sin(t)$ & $0.88$ & $0.92$ \\ \hline
    $\varepsilon(t) = 0.15+0.07\sin(t)$ & $0.08$ & $0.22$ \\ \hline
    $\sigma(t) = 0.12+0.01(\sin(t)+\cos(t))$ &$ 0.12-0.01\sqrt{2}$ &$0.12+0.01\sqrt{2}$\\\hline
    $\Lambda(t) = 0.5+0.06\sin(t)$& $0.44$& $0.56$ \\\hline
    $\mu(t) = 0.07+0.004\cos(t) $ & $0.066$ & $0.074$
 \\\hline
    \end{tabu}
    }
\end{center}
\begin{center}
{\small Table 1: Parameters for simulation 1 of system (\ref{XC}).}
\end{center}

It is important to note that since the initial condition is unchanging this forces two parameters, namely $\Lambda(t)$ and $\mu(t)$, to remain unchanged for these simulation purposes. Moreover, we have that
$$
 \Gamma = \left \{ (x,y,z) \in \mathbb{R}_+^3: x+y+z \le \frac{\overline{\Lambda}}{\underline{\mu}}= 8.484848\right \}
$$
 as the invariant set for the system (\ref{XC}). That is, this system has a unique strong solution taking values in $\Gamma$ per Theorem \ref{RE}. Given these parameters and following Example \ref{exa3.3}, we have
 $$
 	\tilde{R}_0 = \frac{\overline{\beta}\,\overline{\Lambda}}{\underline{\mu}(\underline{\mu}+\underline{\gamma}+
\underline{\varepsilon})}-\frac{\underline{\sigma}^2\overline{\Lambda}^2}{2\underline{\mu}^2(\underline{\mu}+\underline{\gamma}+\underline{\varepsilon})}\ \le 0.7646< 1,\ \ \ \  \underline{\sigma}^2< 0.0121 < 0.0165  = \frac{\underline{\mu}\overline{\beta}}{\overline{\Lambda}}.
 $$
As demonstrated below in Figure 1, the disease will go extinct with exponential rate
$$
	-\alpha \ge (\underline{\mu}+\underline{\gamma}+\underline{\varepsilon})\left( 1- \tilde{R}_0\right) \ge 0.241.
$$

\begin{figure}[H]
\begin{center}
\includegraphics[width=\textwidth]{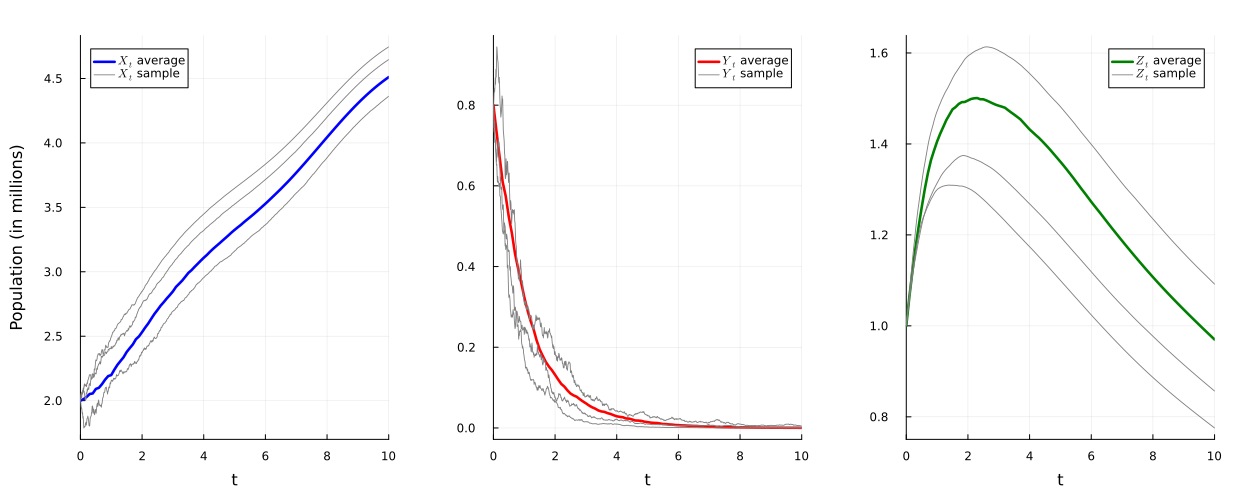}
\caption{Simulation 1 using E-M scheme of system  (\ref{XC}) and displaying the average and three randomly-selected sample paths.}
\end{center}
\end{figure}


We now make only the alteration of a single parameter  in the system (\ref{XC}). Assume that $\sigma(t)$ has the form given in Table 2:

\begin{center}
	\resizebox{9cm}{!}{
    \begin{tabu}{| c  | c | c | }
         \hline
    $f(t)$ &  $\underline{f}$ & $\overline{f}$\\ \hline
  	$\sigma(t) = 0.55+0.003(\sin(t)+\cos(t))$ &$0.55-0.003\sqrt{2}$ &$0.55+0.003\sqrt{2}$\\\hline
    \end{tabu}
    }
\end{center}
\begin{center}
{\small Table 2: Parameters for simulation 2 of system (\ref{XC}).}
\end{center}
\color{black}
This alteration yields
$$
	\underline{\sigma}^2\ge 0.29  > 0.0165 \ge \max \left\{ \frac{\underline{\mu}\overline{\beta}}{\overline{\Lambda}},\ \frac{\overline{\beta}^2}{2(\underline{\mu}+\underline{\gamma}+\underline{\varepsilon})}\right\}.
$$
Thus, we have the scenario in which the disease goes extinct with exponential rate
$$
	-\alpha \ge (\underline{\mu}+\underline{\gamma}+\underline{\varepsilon})- \frac{\overline{\beta}^2}{2\underline{\sigma}^2}\ge 0.993.
$$
Moreover, if we compare Figure 2 to the above Figure 1 we notice the disease appears to go extinct at a faster rate which is as expected given the above results.
\begin{figure}[H]
\begin{center}
\includegraphics[width=\textwidth]{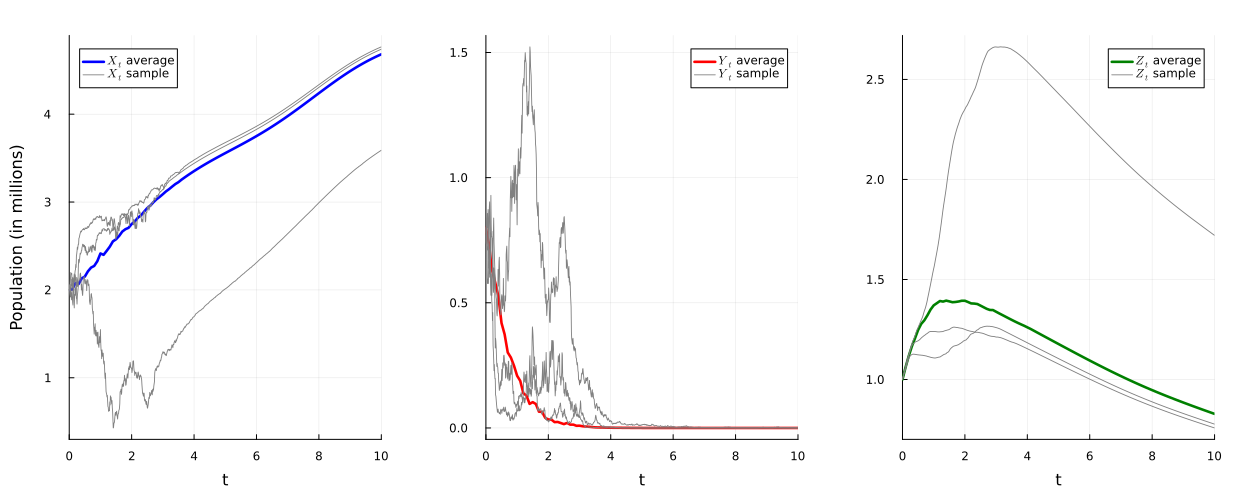}
\caption{Simulation 2 using E-M scheme of system  (\ref{ex1b}) and displaying the average and three randomly-selected sample paths.}
\end{center}
\end{figure}


Now assume that the parameters for the system ($\ref{XC}$) are given in Table 3:
\begin{center}
	\resizebox{9cm}{!}{
    \begin{tabu}{| c  | c | c |}
        \hline
    $f(t)$ &  $\underline{f}$ & $\overline{f}$\\ \hline
    $\beta(t)=0.56+0.01\sin(4t)$ & $0.55$  & $0.57$ \\
    \hline
    $\gamma(t) = 0.25+0. 1\cos(5t)$ & $0.15$ & $0.35$ \\ \hline
    $\sigma(t) = 0.24+0.01(\sin(t)+\cos(t))$ &$0.24-0.01\sqrt{2}$ &$0.24+0.01\sqrt{2}$\\\hline
     \end{tabu}
    }
\end{center}
\begin{center}
{\small Table 3: Parameters for simulation 3 of system (\ref{XC}).}
\end{center}
\color{black}
This modification yields
$$
	\tilde{R}_0 = \frac{\underline{\beta}\underline{\Lambda}}{\overline{\mu}(\overline{\mu}+\overline{\gamma}+\overline{\varepsilon})}
-\frac{\overline{\sigma}^2\overline{\Lambda}^2}{2\underline{\mu}^2(\overline{\mu}+\overline{\gamma}+\overline{\varepsilon})}\ge 1.7 > 1.
$$
In Figure 3 below, we see such a modification yields disease persistence as opposed to disease extinction achieved in the previous two simulations for the system ($\ref{XC}$).
\begin{figure}[H]
\begin{center}
\includegraphics[width=\textwidth]{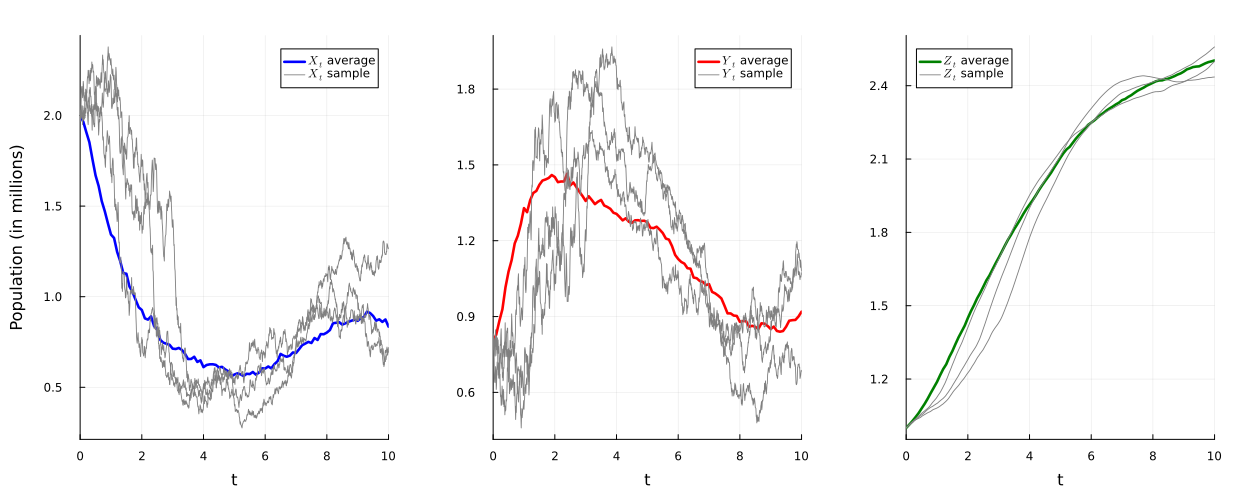}
\caption{Simulation 3 using E-M scheme of system  (\ref{XC}) and displaying the average and three randomly-selected sample paths.}
\end{center}
\end{figure}
\end{simu}
\begin{simu}  We assume the system (\ref{ex34a}) has initial conditions $(3.75, 1.15,1.1)$, where the values are taken to be in millions. Set the parameters in Table $4$:

\begin{center}
	\resizebox{8cm}{!}{
    \begin{tabu}{| c  | c | c |}
        \hline
         $f(t)$ &  $\underline{f}$ & $\overline{f}$\\ \hline
         $M = 2$&  --- & --- \\ \hline
        $\Lambda(t) = 0.15+0.006\sin(t)$& $0.144$& $0.156$\\ \hline
        $\mu(t) = 0.002 + 0.0001\cos(t)$& $0.0019$& $0.0021$\\ \hline
    $\beta(t)=0.18 +0.01\sin(2t)$ & $0.17$  & $0.19$ \\ \hline
    $\gamma_1(t) = 0.15 + 0.004\cos(t) $ & $0.146$ & $0.154$ \\ \hline
    $\gamma_2(t) = 0.12 + 0.02\cos(t) $ & $0.1$ & $0.14$ \\ \hline
     $\gamma_3(t) = 0.12 + 0.04\cos(2t) $ & $0.08$ & $0.16$ \\ \hline
      $\gamma_4(t) = 0.1 + 0.04\sin(4t) $ & $0.06$ & $0.14$ \\ \hline
    $\xi(t) = 1+\ln(1+|\sin(t)|)$ & $1$ & $1+\ln2$ \\ \hline
    $\varphi_i(t) = 0.01 + 0.005\cos(t), i=1,2$ &$0.005$ &$0.015$\\\hline
    $\varphi_3(t) = 1 + 0.25\sin(15t)$ & $0.75$ & $1.25$ \\\hline
    $\sigma_1(t)= 0.015 + 0.01\cos(t) $ & $0.005$& $0.025$\\\hline
    $\sigma_2(t) = 0.012+0.01\sin(t)$& $0.002$& $0.022$\\ \hline
     $h_1(u) = 0.0001$& --- & --- \\ \hline
    $h_2(u) = 0.00025$& --- & --- \\ \hline
    $h_3(u) = 0.0009$& --- & --- \\ \hline
    $g_1(u) = 0.001$& --- & --- \\ \hline
    $g_2(u) = 0.0012$& --- & ---
 \\\hline
    \end{tabu}
    }
\end{center}
\begin{center}
{\small Table 4: Parameters for simulation of system (\ref{ex34a}).}
\end{center}

We have $\lambda_0 = 1.16$ and $\lambda \ge 0.089$ and
$$
\liminf_{t\to\infty}\frac{1}{t}\int_0^tY_sds \ge \frac{0.089}{1.16} \ge 0.076.
$$
The resulting persistence of disease is illustrated below in Figure $4$.
\begin{figure}[H]
\begin{center}
\includegraphics[width=\textwidth]{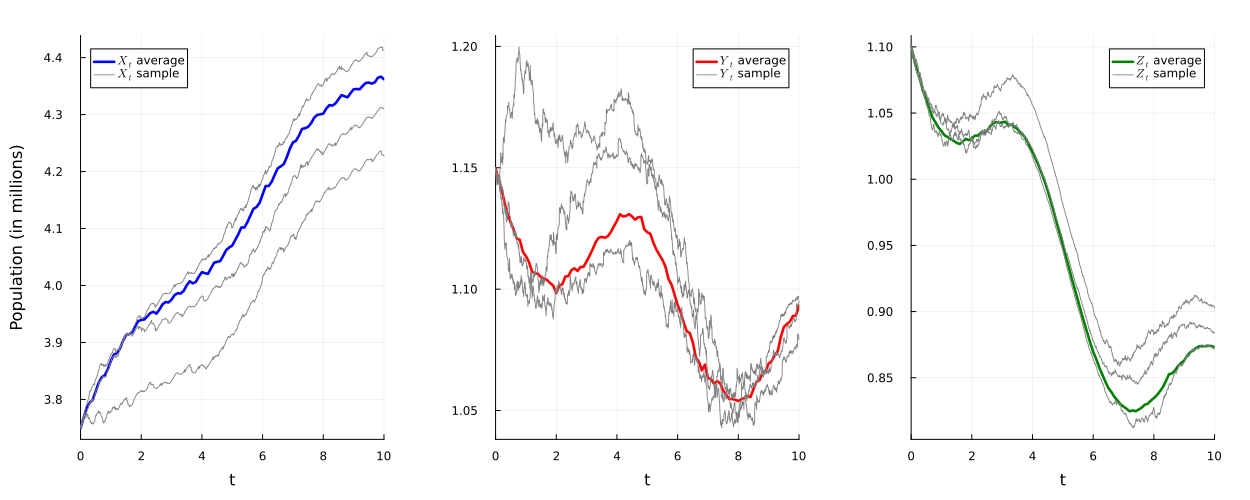}
\caption{Simulation using E-M scheme of system  (\ref{ex34a}) and displaying the average and three randomly-selected sample paths.}
\end{center}
\end{figure}

\end{simu}

\begin{simu} We assume the system (\ref{ex34b}) has initial conditions $(7.27, 1.5,1.11)$, where the values are taken to be in millions. Set the parameters in Table $5$:
\begin{center}
	\resizebox{9cm}{!}{
    \begin{tabu}{| c  | c | c |}
        \hline
         $f(t)$ &  $\underline{f}$ & $\overline{f}$\\ \hline
         $M=1.5$& ---& ---\\ \hline
        $\Lambda(t) = 0.09+0.01\cos(t)$& $0.08$& $0.1$\\ \hline
        $\mu(t) = 0.003+0.001\sin(t)$& $0.002$& $0.004$\\ \hline
    $\beta(t)=0.14+0.005\cos(10t)$ & $0.135$  & $0.145$ \\ \hline
     $\gamma_1(t) = 0.002 + 0.002\cos(25t) $ & $0$ & $0.004$ \\ \hline
    $\gamma_2(t) = 0.35+0.04\cos(15t) $ & $0.31$ & $0.39$ \\ \hline
       $\sigma(t)= 0.3125 + 0.002(\sin(t)+\cos(t)) $ & $0.3125-0.002\sqrt{2}$& $0.3125+0.002\sqrt{2}$\\\hline
    $h_1(u) = 0.0001$& --- & --- \\ \hline
    $h_2(u) = 0.0004$& --- & --- \\ \hline
    $h_3(u) = 0.0009$& --- & --- \\ \hline
    $g_1(u) = 0.001$& --- & --- \\ \hline
    $g_2(u) = 0.007$& --- & --- \\ \hline
    $g_3(u) = 0.005$& --- & ---
 \\\hline
    \end{tabu}
    }
\end{center}
\begin{center}
{\small Table 5: Parameters for simulation of system (\ref{ex34b}).}
\end{center}
The extinction of the disease is illustrated below in Figure $5$. Moreover, as in Example \ref{exa3.4}(b), the disease will go extinct with rate
$-\alpha \ge \underline{\gamma_2} + \underline{\mu} - \overline{\beta} - 2\overline{g_1}\ge 0.165$.
\begin{figure}[H]
\begin{center}
\includegraphics[width=\textwidth]{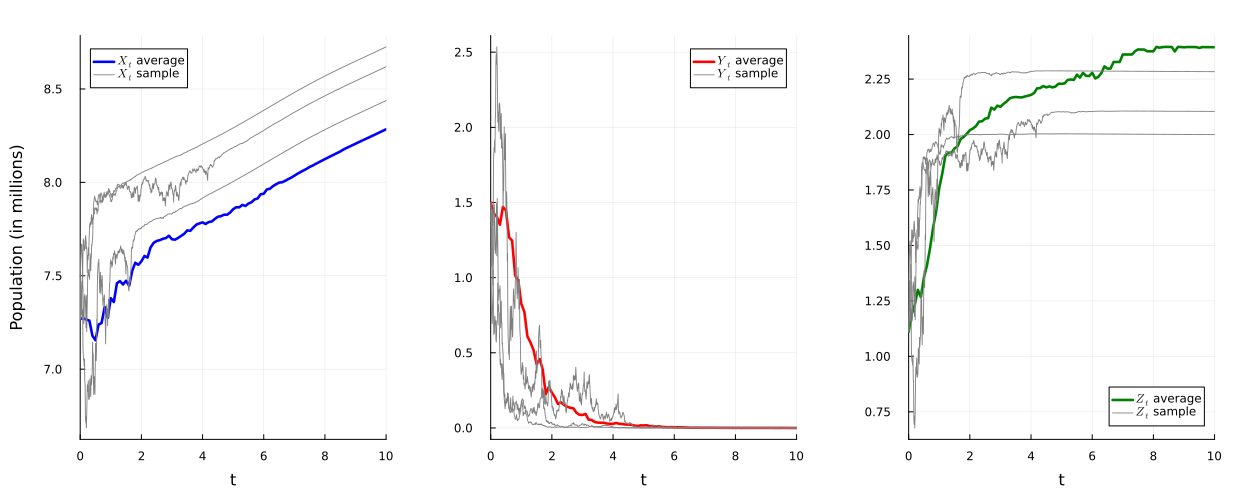}
\caption{Simulation using E-M scheme of system  (\ref{ex34b}) and displaying the average and three randomly-selected sample paths.}
\end{center}
\end{figure}
\end{simu}
\begin{simu} We assume that the system (\ref{p1_ex1}) in Example \ref{exa2.5}(a) has initial values $(X_0,Y_0,Z_0)=(0.8,0.19, 0.01)$ and set the parameters in Table 6:

\begin{center}
	\resizebox{9cm}{!}
	{
    \begin{tabu}{| c  | c | c |}
        \hline
    $f(t)$ &  $\underline{f}$ & $\overline{f}$\\ \hline
    $\beta(t)=0.3+0.1\sin(4t)$ & $0.2$  & $0.4$ \\
    \hline
    $\gamma(t) = 0.8 + 0.04\cos(7t) $ & $0.76$ & $0.84$ \\ \hline
    $\xi(t) = 1+\frac{t}{1+t}$ & $1$ & $2$ \\ \hline
    $\varphi_i(t) = 0.01 + 0.005\cos(t), i=1,2$ &$0.005$ &$0.015$\\\hline
    $\varphi_3(t) = 1+0.5\sin(15t)$ & $0.5$ & $1.5$ \\\hline
    $\sigma_1(t)= 0.5 + 0.01\cos(7t) $ & $0.49$& $0.51$\\\hline
    $\sigma_2(t) = 0.4+0.01\sin(7t)$& $0.39$& $0.41$\\ \hline
      $h_1(u) = 0.01$& --- & --- \\ \hline
    $h_2(u) = 0.025$& --- & --- \\ \hline
    $g_1(u) = 0.1$& --- & --- \\ \hline
    $g_2(u) = 0.12$& --- & ---
 \\\hline
    \end{tabu}
    }
\end{center}
\begin{center}
{\small Table 6: Parameters for simulation of system (\ref{p1_ex1}).}
\end{center}
In Figure 6 below, it is illustrated that the extinction of the disease occurs at an exponential rate. In accordance with Example \ref{exa2.5}(a), the disease will extinct with exponential rate
$$-\alpha \ge \underline{\gamma} - \overline{\beta} - 2\overline{g_1} \ge 0.16.$$
\begin{figure}[H]\label{fig1}
\begin{center}
\includegraphics[width=\textwidth]{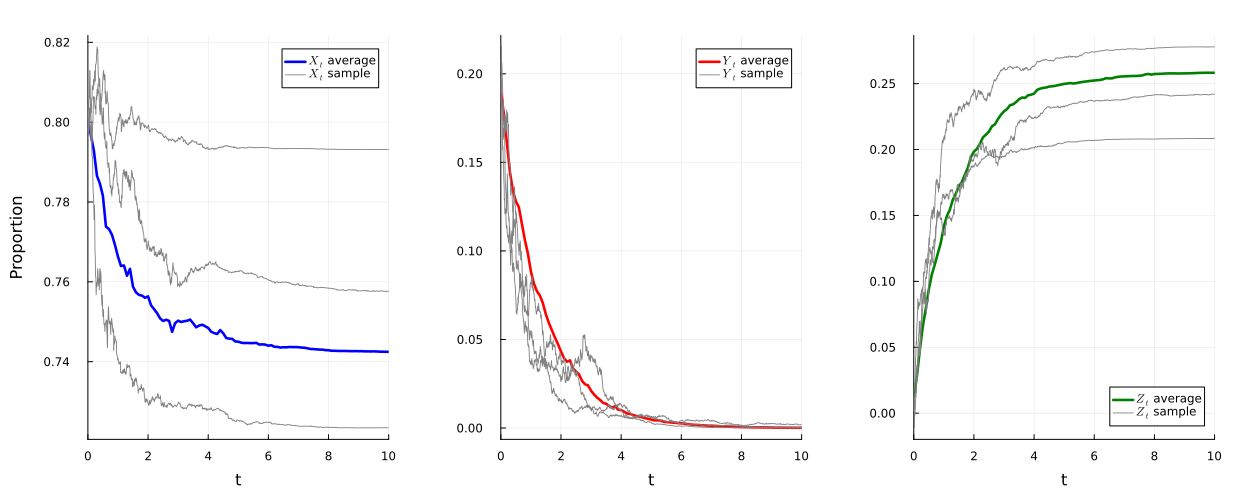}
\caption{Simulation using E-M scheme of system  (\ref{p1_ex1}) and displaying the average and three randomly-selected sample paths.}
\end{center}
\end{figure}
\end{simu}
\begin{simu} We assume that the system (\ref{ex1b}) in Example \ref{exa2.5}(b) has initial values $(X_0,Y_0,Z_0)=(0.85, 0.1, 0.05)$ and set the parameters in Table 7:
\begin{center}
	\resizebox{9cm}{!}{
    \begin{tabu}{| c  | c | c |}
       \hline
    $f(t)$ &  $\underline{f}$ & $\overline{f}$\\ \hline
    $\beta(t)=0.17+ 0.01\cos(20t)$ & $0.16$  & $0.18$ \\
    \hline
    $\gamma_1(t) = 0.12+ 0.01\cos(t) $ & $0.11$ & $0.13$ \\ \hline
    $\gamma_2(t) = 0.56+0.01\sin(t) $ & $0.55$ & $0.57$ \\ \hline
    $\sigma_1(t)=\sigma_2(t)=\sigma_3(t)=0.141+0.02(\sin(t) + \cos(t)) $& $0.141-0.02\sqrt{2}  $& $0.141+0.02\sqrt{2}$ \\ \hline
    $h_1(u) = 0.019$& --- & --- \\ \hline
    $h_2(u) = 0.018$& --- & --- \\ \hline
    $h_3(u) = 0.018$& --- & --- \\ \hline
    $g_1(u) = 0.11$& --- & --- \\ \hline
    $g_2(u) = 0.01$& --- & --- \\ \hline
    $g_3(u) = 0.01$& --- & ---
 \\\hline
    \end{tabu}
    }
\end{center}
\begin{center}
{\small Table 7: Parameters for simulation of system (\ref{ex1b}).}
\end{center}

We achieve results which illustrate persistence of the disease, as is displayed in Figure 7 below.  Furthermore, we have $\lambda_0 = 0.55 $, $\lambda= 0.21  $
and $$ \liminf_{t\to \infty } \frac{1}{t} \int_0^t Y_sds \ge \frac{0.21 }{0.55} \ge 0.38. $$

\begin{figure}[H]
\begin{center}
\includegraphics[width=\textwidth]{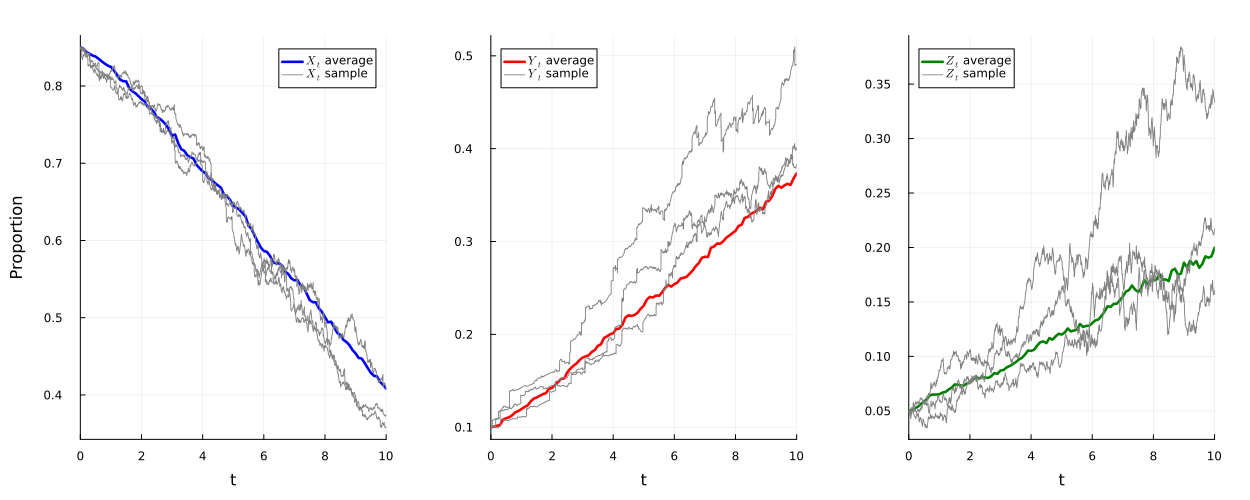}
\caption{Simulation using E-M scheme of system  (\ref{ex1b}) and displaying the average and three randomly-selected sample paths.}
\end{center}
\end{figure}

\end{simu}

\section{Conclusion}
In this paper, we propose and investigate the USSIR model given by the system (\ref{p1_eq2}). We have presented two forms of the novel model -- one for population  numbers and the other for population proportions. For both forms of the model, we have given results on the extinction and persistence of  diseases; moreover, we have shown that these results still hold with time-dependent, nonlinear parameters and multiple L\'evy noise sources. Notably, we give examples and simulations that agree with the theoretical results and illustrate the impact that noise has on a given SIR model system. Moreover, the ability to allow time-dependency and multiple noises coincides with real world occurrences of infectious disease spread due to environmental noises or time-dependent events such as temperature, climates, seasons, and so forth. Our examples are intended to have an emphasis on real-world relevance; however, we do surmise the inclusion of artificiality in examples illustrates the flexibility of our model.

There is much work to follow this introduction to the USSIR model. An initial follow-up problem is parameter estimation of the USSIR model including the presence of periodic parameters. Additionally, given the measurement of real-world data is distorted by noise and may contain unknowns, it is important to consider the USSIR model with filtering to overcome these difficulties and inaccuracies. Both of these problems will be considered in future work with applications to real-world data. Additionally, the methods utilized here certainly would be applicable to models of dimension larger than $3$ which may be explored in further work.

\section*{\large \center Acknowledgements}
We would like to thank the anonymous reviewers and the Editor
for suggestions and comments leading to a significant improvement
of the paper. Also, we wish to thank Professor Yang Lu for his constructive and insightful comments on an earlier version of the paper. This work  was partially supported by the Natural Sciences and Engineering Research Council of Canada (No.
4394-2018).

\end{document}